\documentclass[12pt, twoside, reqno]{amsart}
\calclayout
\usepackage[utf8]{inputenc}
\usepackage{amsmath,amssymb,amscd,amsxtra,latexsym,bbm,graphicx,epsfig,epic,
eepic,oldgerm,bm,psfrag,amsthm,mathrsfs,bm,bbm, verbatim}
\usepackage{caption}
\usepackage{subcaption}
\usepackage{tikz-cd}
\usepackage{enumitem}
\usepackage[pdftex]{hyperref}
\usetikzlibrary{shapes}

\tikzset{cross/.style={
    cross out, draw, solid, thin, 
    minimum size=2*(#1-\pgflinewidth), 
    inner sep=0pt, outer sep=0pt
  },
  cross/.default={3}
}

\input xy
\xyoption{all}
\CompileMatrices
    
\theoremstyle{plain}
\newtheorem{theorem}{Theorem}[section]
\newtheorem*{theorem*}{Theorem}
\newtheorem{lemma}[theorem]{Lemma}
\newtheorem{proposition}[theorem]{Proposition}
\newtheorem{corollary}[theorem]{Corollary}
\newtheorem{definition}[theorem]{Definition}
\newtheorem*{definition*}{Definition}

\newtheorem{question}[theorem]{Question}

\theoremstyle{definition}
\newtheorem*{remark*}{Remark}
\newtheorem*{remarks*}{Remarks}
\newtheorem{remark}[theorem]{Remark}

\newtheorem{example}[theorem]{Example}

\newtheorem*{example*}{Example}
\newtheorem*{examples*}{Examples}

\theoremstyle{plain}
\newtheorem{TheoremA}{Theorem}

\newcommand{\proofend}{\hspace*{\fill} $\Box$\\}

\newcommand{\ign}[1]{}

\def\1{\:\!}
\def\2{\;\!}

\def\im{\operatorname {im}}

\def\Diffc0{\operatorname{Diff^c_0}}
\def\Symp{\operatorname{Symp}}

\def\Sympc0{\operatorname{Symp^c_0}}

\def\Int{\operatorname{int}}
\def\Ham{\operatorname{Ham}}

\def\Crit{\operatorname{Crit}}

\def\Flux{\operatorname{Flux}}

\def\GL{\operatorname{GL}}

\def\Span{\operatorname{span}}

\def\IA{\operatorname{IA}}

\def\CP{\operatorname{CP}}

\def\ce{{\mathcal E}}

\def\cj{{\mathcal J}}

\def\cl{{\mathcal L}}

\def\cu{{\mathcal U}}

\def\fm{{\mathfrak{m}}}

\def\ob{\overline{b}}

\def\C{\mathbb{C}}

\def\N{\mathbb{N}}

\def\R{\mathbb{R}}

\def\Z{\mathbb{Z}}

\def\RP{\R P}
\def\CP{\C P}

\def\pp{\partial}

\def\ddt0{\left. \frac{d}{dt} \right\vert_{t=0}}
\def\dds0{\left. \frac{d}{ds} \right\vert_{s=0}}
\def\ddt{\frac{d}{dt} }
\def\dds{\frac{d}{ds} }

\def\^t{^{\times}}

\def\ni{\noindent}

\def\.{\mskip1mu}
\def\?{\mskip-1mu}

\def\id{\operatorname{id}}

\def\proof{\noindent {\it Proof. \;}}
\newcommand{\proofof}[1]{\ni {\it Proof of #1. }}

\begin{document}

\title[]{Local exotic tori}

\author{Jo\'e Brendel}  

\address{Jo\'e Brendel,
ETH Zürich,
D-MATH, 
Rämistrasse 101,
8092 Zürich, 
Switzerland }
\email{joe.brendel@math.ethz.ch}

\date{\today}

\begin{abstract} 
For a broad class of symplectic manifolds of dimension at least six, we find the following new phenomenon: there exist local exotic Lagrangian tori.

More specifically, let~$X$ be a geometrically bounded symplectic manifold of dimension at least six. We show that every open subset of~$X$ contains infinitely many Lagrangian tori which are distinct up to symplectomorphisms of~$X$ while being Lagrangian isotopic and having the same classical invariants. The proof relies on a locality property of the displacement energy germ, which allows us to compute it in a Darboux chart. 

Since these tori are not monotone, bubbling may occur and the count of Maslov index two $J$-holomorphic disks does not yield an invariant. 
\end{abstract}

\maketitle

\section{Introduction}

This paper is motivated by the classification question of Lagrangian submanifolds up to symplectomorphisms of the ambient symplectic manifold. The full classification is open in most cases, for example for Lagrangian tori in~$\R^4$. In this paper we show that infinitely many Lagrangian tori, which cannot be mapped to one another by an ambient symplectomorphism, occur in every open set of every geometrically bounded symplectic manifold of dimension at least six.

\subsection{Lagrangian tori in symplectic vector spaces} Equip~$\R^2 = \C$ with the standard area form~$\omega_0 = d\lambda_0$. In~$\R^2$, the only closed Lagrangians are closed embedded curves and Hamiltonian isotopies are just smooth area-preserving isotopies with compact support. Every closed embedded curve~$L$ is isotopic to a circle~$S^1(a) = \{\pi \vert z \vert^2 = a\}$ by an area-preserving isotopy which has compact support, where~$a = \int_L \lambda_0$ is the area enclosed by~$L$. This settles the classification question of Lagrangian tori in~$(\R^2, \omega_0)$. 

Let~$\R^{2n} = \C^n$ be equipped with the product symplectic form~$\omega_0 \oplus \ldots \oplus \omega_0 = d\lambda$. Then we obtain a family of Lagrangian tori, called \emph{product tori} (or \emph{split tori}), by setting
\begin{equation}
	\label{eq:proddtorii}
	T(a_1,\ldots,a_n) = S^1(a_1) \times \ldots \times S^1(a_n), \quad
	a_i > 0. 
\end{equation}
Inspired by the case~$n = 1$, it is natural to ask if every Lagrangian torus in~$\C^n$ is Hamiltonian isotopic to a product torus. 

In~$\C^2$ this question was answered negatively by Chekanov~\cite{Che96}, who constructed the so-called \emph{Chekanov torus}, which cannot be mapped by a symplectomorphism to a product torus. A Lagrangian torus in a symplectic vector space~$\R^{2n} = \C^n$ which cannot be mapped by a symplectomorphism to a product torus is called \emph{exotic}. An alternate construction of the Chekanov torus, which heavily influenced the present paper, was given by Eliashberg--Polterovich~\cite{EliPol97}. See also~\cite{Aur07, Bre22, CheSch10} for different points of view on the same construction. Despite significant progress made by Dimitroglou--Rizell~\cite{Riz19} towards a classification, it is unknown as of yet if there are other tori in~$\R^4$ besides product and Chekanov tori. 

In~$\C^3$, on the other hand, it is known that there are infinitely many Lagrangian tori which are exotic and pairwise distinct up to symplectomorphisms. This was proved by Auroux~\cite{Aur15}. The Chekanov torus as well as Auroux's examples are \emph{monotone} (meaning their area and Maslov classes are positively proportional) and appear in one-parameter families. For every~$k \in \Z_{\geqslant 2}$, set
\begin{equation}
	\label{eq:uk}
	U_k = \left\{(a_1,a_2) \in \R_{>0}^2 \, \left\vert \, \frac{a_1}{k+1} \leqslant a_2 < \frac{a_1}{k} \right. \right\} \subset \R^2.
\end{equation}
In~\S\S\ref{ssec:famsympred}-\ref{ssec:torivd}, we construct Lagrangian tori,
\begin{equation}
	\label{eq:thmAeq}
	\Upsilon_k(a) \subset \C^3, \quad
	a \in U_k.
\end{equation}
The parameter~$a=(a_1,a_2)$ has the following geometric significance. There is a basis~$e_1,e_2,e_3 \in H_1(\Upsilon_k(a))$ of Maslov two classes with
\begin{equation}
	\label{eq:upsilonarea}
	\int_{e_1}\lambda\vert_{\Upsilon_k(a)} = a_1 - ka_2, \quad
	\int_{e_2} \lambda\vert_{\Upsilon_k(a)} = \int_{e_3} \lambda\vert_{\Upsilon_k(a)} = a_2. 
\end{equation}
This means that~$a$ determines the area class as defined in~\eqref{eq:softinvvv}. In particular, the torus~$\Upsilon_k(a)$ is monotone if and only if~$a_1 = (k+1)a_2$. In~\S\ref{ssec:proofthmA}, we prove the following.

\begin{TheoremA}
\label{thm:A}
For every~$k \in \Z_{\geqslant 2}$, the torus~$\Upsilon_k(a)$ is exotic. Let~$k \neq k'$ and~$a \in U_k, a' \in U_{k'}$, then~$\Upsilon_k(a)$ cannot be mapped to~$\Upsilon_{k'}(a')$ by a symplectomorphism of~$\C^3$.
\end{TheoremA}

Our construction, as well as our method of distinguishing Lagrangians differs from~\cite{Aur15}. Due to the occurence of non-monotone exotic tori, Theorem~\ref{thm:A} strenghtens Auroux's result.

\begin{remark}
We do not know for which~$a,a' \in U_k$, the tori~$\Upsilon_k(a)$ and~$\Upsilon_k(a')$ are equivalent up to a symplectomorphism of~$\C^3$. Note that even the classification of product tori up to ambient symplectomorphisms is a subtle question, see~\cite[Theorem A]{Che96}.
\end{remark}

\begin{remark}
\label{rk:vianna}
The monotone members of~\eqref{eq:thmAeq} were independently constructed and studied in unpublished work by Dmitry Tonkonog and Renato Vianna; Jeff Hicks was also aware of the construction. Furthermore, R. Vianna informed the author that it is likely, based on computing the superpotential, that the monotone subfamily of the~$\Upsilon_k$ corresponds to Auroux's examples~\cite{Aur15}.
\end{remark}

\subsection{Lagrangian tori in general symplectic manifolds}
We consider the question of exotic tori in more general ambient spaces. In sharp contrast to what is known in the case of~$\C^2$, Vianna~\cite{Via17, Via16} proved that there are infinitely many monotone Lagrangian tori in~$\CP^2$ and other Del Pezzo surfaces. We will discuss this in detail in~\S\ref{ssec:atfibres}. This result was conjectured earlier by Galkin--Usnich~\cite{GalUsn10} and worked out independently by Galkin--Mikhalkin~\cite{GalMik22}. Earlier work producing finitely many exotic Lagrangian tori in different ambient spaces includes~\cite{AlbFra08,BirCor09, CheSch10,EntPol09,FOOO12,Gad13,OakUsh16,Wu15}. We refer to the introduction of~\cite{OakUsh16} for an overview.

The main aim of this paper is to construct \emph{exotic tori} in general (geometrically bounded) symplectic manifolds of dimension greater or equal to six. The notion of \emph{exoticity} only makes sense in spaces in which there is a natural family of Lagrangian tori, from which the \emph{exotic} examples are distinguished. For example, in~$\C^n$, the standard tori are product tori; in more general toric manifolds, standard tori are toric fibres. 

For a general symplectic manifold~$(X,\omega)$, let us clarify the question we are interested in.

\begin{definition}
We say that two Lagrangians~$L,L' \subset (X,\omega)$ are \emph{equivalent} if there is~$\phi \in \Ham(X,\omega)$ with~$\phi(L) = L'$. This is denoted by~$L \cong L'$. We say that~$L, L'$ are \emph{inequivalent} or \emph{distinct} if there is no~$\phi \in \Symp(X,\omega)$ with~$\phi(L) = L'$ and denote this by~$L \ncong L'$. 
\end{definition}

The mix of Hamiltonian diffeomorphisms and symplectomorphisms in the definition is justified by the fact that in all the examples we treat, Lagrangians are distinguished in the strongest sense (up to symplectomorphisms) or shown to be equivalent in the strongest sense (up to Hamiltonian diffeomorphisms). We are thus looking for sets of \emph{distinct} Lagrangian tori. 

Lagrangian tori may be distinct because they have different \emph{soft invariants}. By soft invariants, we mean Maslov and area class,
\begin{equation}
	\label{eq:softinvvv}
	m_L \colon H_2(X,L) \rightarrow \Z, \quad
	\sigma_L \colon H_2(X,L) \rightarrow \R.
\end{equation}
We have~$\phi^* m_{\phi(L)} = m_L$ and~$\phi^* \sigma_{\phi(L)} = \sigma_L$ for all~$\phi \in \Symp(X,\omega)$. For example if~$\im \sigma_L \neq \im \sigma_{L'} \subset \R$, then~$L \ncong L'$. In light of this, let us make the following definition. 

\begin{definition}
\label{def:weakequiv}
We say that~$L,L' \subset X$ are \emph{area-equivalent (in $X$)} if there is a Lagrangian isotopy with endpoints~$L$ and~$L'$ preserving the area class in the sense that~$\Phi^* \sigma_{L'} = \sigma_L$, where~$\Phi \colon H_2(X,L) \rightarrow H_2(X,L')$ is the isomorphism induced by the Lagrangian isotopy. We write~$L \sim L'$ in that case.
\end{definition}

The Maslov class is preserved by Lagrangian isotopies, meaning that both soft invariants of~$L$ and~$L'$ agree if~$L$ and~$L'$ are area-equivalent. Both the area and the Maslov class are preserved by Hamiltonian isotopies (and symplectomorphisms more generally), meaning that~$L \cong L'$ implies~$L \sim L'$. The converse gives rise to the main question addressed in this paper. 

\begin{question}
\label{q:main}
For a given symplectic manifold~$X$, are there Lagrangians~$L,L' \subset X$ such that~$L \sim L'$, but~$L \ncong L'$? 
\end{question}

All of the exotic tori from~\cite{Aur15, Che96, CheSch10, Via17, Via16} discussed so far (and all other examples of exotic tori known to the author) yield examples answering Question~\ref{q:main} positively. For example the Chekanov torus in~$\C^2$ is area-equivalent to a monotone product torus~$T(a,a)$ for a suitable~$a > 0$, depending on the size of the chosen Chekanov torus. Similarly, the tori~$\Upsilon_k(a)$ appearing in Theorem~\ref{thm:A} are area-equivalent for different~$k$. More precisely, we prove that~$\Upsilon_k(a_1,a_2) \sim T(a_1 -k a_2, a_2, a_2)$, where~$T(\cdot)$ denotes a product torus as in~\eqref{eq:proddtorii}. See also Remark~\ref{rk:areaequivalent}.

We prove the following Theorem in~\S\ref{ssec:proofthmB}. See~\cite[Chapter X, Definition 2.2.1]{AudLaf94} for a definition of geometrically bounded symplectic manifolds.

\begin{TheoremA}
\label{thm:B}
Let~$(X^6,\omega)$ be a geometrically bounded symplectic manifold and~$\psi \colon B^6(R) \rightarrow X$ a Darboux chart. Then there is~$\varepsilon > 0$ such that 
\begin{equation}
	\psi(\Upsilon_k(a)) \ncong \psi(\Upsilon_{k'}(a'))
\end{equation}
for all~$k \neq k'$ and~$a = (a_1,a_2) \in U_k, a'=(a_1',a_2') \in U_{k'}$ with~$a_1,a_1' < \varepsilon$.
\end{TheoremA}

Again, we have~$\psi(\Upsilon_k(a)) \sim \psi(T(a_1 - k a_2, a_2, a_2))$ in~$X$, meaning that we find tori answering Question~\ref{q:main} in the affirmative.

The condition~$a_1 < \varepsilon$ means that the tori to which Theorem~\ref{thm:B} applies have a sufficiently small area class. The threshold~$\varepsilon > 0$ is discussed in Remark~\ref{rk:epsilon}. A significant drawback of this construction is the following. For fixed~$\varepsilon > 0$, the area class of tori to which the statement of Theorem~\ref{thm:B} applies become smaller as~$k$ increases. Indeed, this follows from~\eqref{eq:upsilonarea}, which yields the area class of~$\Upsilon_k(a)$ in terms of~$a \in U_k$, togther with~\eqref{eq:uk}. This means that in a fixed area class, Theorem~\ref{thm:B} produces only finitely many inequivalent Lagrangian tori in~$(X,\omega)$.\smallskip

To remedy this drawback, we construct another family of exotic tori in~$\C^3$ by lifting Vianna's tori in~$\CP^2$ to tori in~$\C^3$ via the Hopf reduction map. For every Markov triple~$\fm$, meaning every triple~$\fm = (\alpha,\beta,\gamma)$ of natural numbers satisfying the Diophantine equation
\begin{equation}
	\label{eq:markoveq}
	\alpha^2 + \beta^2 + \gamma^2 = 3\alpha\beta\gamma,
\end{equation}
there is a monotone Lagrangian torus~$T_{\fm} \subset \CP^2$ appearing as the monotone fibre of a certain almost toric fibration of~$\CP^2$ with~$T_{\fm} \cong T_{\fm'}$ if and only if~$\fm = \fm'$, see~\S\S\ref{ssec:cp2mutations}-\ref{ssec:atfibres} or~\cite{GalMik22, Via16}. In~\S\ref{ssec:viannalift}, we show that the lifted tori~$\Theta_{\fm}(a) = p_{a}^{-1}(T_{\fm}) \subset \C^3$ by the reduction map~$p_a \colon \C^3 \supset S^5(a)\rightarrow (\CP^2, a\omega_{\CP^2})$ have the same property. Note that these tori are monotone, as well. 

\begin{proposition}
\label{prop:viannalift}
Let~$\fm, \fm'$ be Markov triples, then $\Theta_{\fm}(a) \cong \Theta_{\fm'}(a)$ if and only if~$\fm = \fm'$.
\end{proposition}

This result was expected and may have been known to experts. What is interesting in the context of the present paper is that this property persists under Darboux embeddings into geometrically bounded symplectic manifolds, provided~$a > 0$ is small enough.

\begin{TheoremA}
\label{thm:C}
Let~$(X^6,\omega)$ be a geometrically bounded symplectic manifold and~$\psi \colon B^6(R) \rightarrow X$ a Darboux chart. Then there is~$\varepsilon > 0$ such that 
\begin{equation}
	\psi(\Theta_{\fm}(a)) \ncong \psi(\Theta_{\fm'}(a)), \quad
	\psi(\Theta_{\fm}(a)) \sim \psi(\Theta_{\fm'}(a))
\end{equation}
for all~$a < \varepsilon$ and~$\fm \neq \fm'$.
\end{TheoremA}

See~\S\ref{ssec:proofthmC} for a proof. Since the set of Markov triples is infinite, Theorem~\ref{thm:C} shows that~$(X,\omega)$ contains infinitely many distinct Lagrangian tori in a (small enough) fixed area class. Thus Theorem~\ref{thm:C} remedies the drawback of Theorem~\ref{thm:B} pointed out above. The idea used to upgrade Proposition~\ref{prop:viannalift} to Theorem~\ref{thm:C} is the same as the one used to upgrade Theorem~\ref{thm:A} to Theorem~\ref{thm:B}, and will be discussed in detail in~\S\ref{ssec:disttori}. Again, we refer to Remark~\ref{rk:epsilon} for a discussion of the size~$\varepsilon > 0$ of tori to which Theorem~\ref{thm:C} applies.

\begin{remark}
\label{rk:nonmonotonelifts}
Note that, although we will not consider these in detail, the lifted Vianna tori~$\Theta_{\fm}$ can be extended to natural two-parameter families, by lifting non-monotone Vianna tori from~$\CP^2$ to~$\C^3$. These were for example considered in a paper by Shelukhin--Tonkonog--Vianna~\cite{SheTonVia19}, see also Example~\ref{ex:nonmonvianna}. Although these two-parameter families can be distinguished by their displacement energy germ by the methods used in this paper, we restrict our attention to the one-parameter families~$\{\Theta_{\fm}(a)\}_{a > 0}$ coming from lifts of monotone Vianna tori to simplify the exposition.
\end{remark}

\begin{remark}
The families~$\Upsilon_k$ and~$\Theta_{\fm}$ are disjoint up to one exception in the sense that~$\Upsilon_k \ncong \Theta_{\fm}$ for all~$k \geqslant 3$ and all Markov triples~$\fm \neq (2,1,1)$. The exception is~$k = 2$ and~$\fm = (2,1,1)$, in which case~$\Upsilon_2(3a,a) \cong \Theta_{(2,1,1)}(3a)$ for all~$a > 0$. The same statements hold for the embeddings of the tori into geometrically bounded symplectic manifolds by Darboux charts. The distinction of the tori follows from the fact that they have inequivalent displacement energy germs, see~\S\ref{ssec:proofthmA} and~\S\ref{ssec:viannalift} for~$\Upsilon_k, \Theta_{\fm} \subset \C^3$ and~\S\ref{ssec:proofthmB} and~\S\ref{ssec:proofthmC} for their embeddings into tame symplectic manifolds. 
\end{remark}

\subsection{Higher dimensions}
\label{ssec:higherdim}

For simplicity, we have restricted our attention to dimension six in the above discussion. Since our methods behave well when taking products, one obtains the statements of Theorems~\ref{thm:A},~\ref{thm:B} and~\ref{thm:C} in higher dimensions by considering tori of the form~$L \times T(c_1,\ldots,c_m) \subset \C^{m + 3}$, where~$L \subset \C^3$ is a~$\Upsilon_k$- or~$\Theta_{\fm}$-torus. The so-obtained tori can be distinguished from one another, provided the~$c_i$ are large enough for the displacement energy germ of the~$L$-factor to be visible in the displacement energy germ of the product, see Remark~\ref{rk:asmall} for details. Furthermore, similar methods can be used to distinguish products of exotic tori of either~$\Upsilon_k$- or~$\Theta_{\fm}$-type. Using products, we prove the following.

\begin{TheoremA}[Local exotic tori]
\label{thm:D}
Let~$X$ be a geometrically bounded symplectic manifold of dimension at least six and~$U \subset X$ an open subset. Then~$U$ contains infinitely many area-equivalent tori which are distinct up to symplectomorphisms of~$X$.
\end{TheoremA}

These local exotic tori naturally appear in~$(n-1)$-dimensional families, where~$\dim X = 2n$, which seems to be a general feature of exotic tori, see Remark~\ref{rk:codimoneee} and Question~\ref{q:2}. Beyond Theorem~\ref{thm:D}, we do not give a comprehensive treatment of tori obtained as products from known tori. All proofs in that direction follow the same pattern as the proof of Theorem~\ref{thm:D} given in~\S\ref{ssec:proofthmD}.

\subsection{Key ideas: construction}
\label{ssec:introkey}

The idea for the construction of the tori~$\Upsilon_k$ is inspired by the Eliashberg--Polterovich~\cite{EliPol97} construction of the Chekanov torus in~$\C^2$ and its subsequent interpretations by Auroux~\cite{Aur07, Aur09} and Chekanov--Schlenk~\cite{CheSch10}. Let us review this construction of the Chekanov torus in~$\C^2$ using the language of toric geometry and symplectic reduction and illustrate how it inspires the construction of the~$\Upsilon_k$-tori. For more details on the point of view discussed here, we refer to the expository paper~\cite{Bre22}.

\begin{figure}
		\begin{tikzpicture}[scale=0.85]	
			\fill[black!10] (0,0)--(3,5)--(14,5)--(11,0)--(0,0);
			\draw [thick,dotted] (0,0)--(3,5)--(14,5)--(11,0);
			\draw [black,->] (0,0)--(0.75,1.25);
			\draw [black,->] (0,0)--(0,5);
			\draw [black,->] (0,0)--(12,0);
			\draw [black!50, ->] (0,0)--(3,2);
			
			\draw [thick, black] (10,0)--(1.5,2.5);
			\draw [thick, black] (8,0)--(2.25,3.75);
			
			\fill[thick, black] (10,0)  circle[radius=1.5pt];	
			\fill[thick, black] (1.5,2.5)  circle[radius=1.5pt];	
			\fill[thick, black] (8,0)  circle[radius=1.5pt];	
			\fill[thick, black] (2.25,3.75)  circle[radius=1.5pt];	
			
			\node at (10,-0.4){$a_1'$};
			\node at (8,-0.4){$a_1$};
			
			\node at (4.2,3.25){$\ell_2(a_1)$};
			\node at (8.95,0.75){$\ell_4(a_1')$};
			
			\node at (1.65,1.25){$(0,1,1)$};
			
			\node at (12,-0.4){$x$};
			\node at (-0.35,4.9){$z$};
			\node[black!65] at (3,1.6){$y$};
			
		\end{tikzpicture}
	\caption{Examples of level sets $\ell_k(a_1)$ on the preimage of which symplectic reduction is performed.}
	\label{fig:1}
\end{figure}

Let~$\ell = \{(x,x)\} \subset \R^2_{\geqslant 0}$ be the diagonal in the moment map image~$\R^2_{\geqslant 0}$ of~$\C^2$ under the standard moment map~$\mu(z_1,z_2) = (\pi \vert z_1 \vert^2, \pi \vert z_2 \vert^2)$. Note that the regular fibres of this moment map are precisely the product tori~$T(a_1,a_2) \subset \C^2$. To the ray~$\ell$, we can naturally associate a symplectic quotient space\footnote{In general, every rational affine subspace~$V \subset \R^n$ intersected with a toric moment polytope~$\Delta \subset \R^n$ gives rise to a (possibly singular) symplectic reduction by the action of the orthogonal torus~$T_V = \exp(V^{\perp}) \subset T^n$ and has as reduced space a (possibly singular) toric symplectic manifold. We refer to~\cite{Bre23} for details on \emph{toric reduction}.} by viewing~$\ell$ as the projection of the level set~$\{ H = 0\}$ of the Hamiltonian~$H = \pi \vert z_1 \vert^2 - \pi \vert z_2 \vert^2$ generating the Hamiltonian~$S^1$-action \emph{orthogonal} to~$\ell$. Note that the level set~$\{H = 0\}$ contains the fixed point~$0 \in \C^2$, which we remove to obtain the reduced space~$(\{H = 0\} \setminus \{0\}) / S^1$ symplectomorphic to the punctured standard plane~$(\C^{\times}, \omega_0\vert_{\C^{\times}})$. On one hand, for~$a>0$ the product torus~$T(a,a) \subset \{H=0\}$ projects to the circle~$S^1(a) \subset \C^{\times}$. On the other hand one can lift closed embedded curves from the reduced space to obtain tori in~$\C^2$. The Chekanov torus~$T_{\rm Ch}^2(a)$ is obtained by lifting a curve in~$\C^{\times}$ bounding a disk of symplectic area~$a > 0$ in~$\C^{\times}$. Note that such circles are not isotopic to standard circles~$S^1(a)$ in~$\C^{\times}$, since the latter enclose the puncture and the former do not. 

The tori~$\Upsilon_k(a)$ are similarly defined by choosing suitable families of rational segments
\begin{equation}
	\ell_k(a_1) = \{ (a_1,0,0) + t(-k,1,1) \,\vert\, t \in \R \} \cap \R^3_{\geqslant 0}, \quad
	k \in \Z_{\geqslant 2}, a_1 > 0
\end{equation}
in the moment polytope~$\R^3_{\geqslant 0}$ of~$\C^3$ and performing~$T^2$-reduction on their preimages. See Figure~\ref{fig:1} for an illustration. The reduced space is a sphere of area~$\frac{a_1}{k}$ with a puncture and an orbifold point\footnote{The puncture comes from removing a set of points with $S^1$-stabilizer and the orbifold point to a set of points with $\Z_k$-stabilizer.} of order~$k$ and the torus~$\Upsilon_k(a_1,a_2)$ is defined as the lift of a closed embedded curve $C(a_2)$ bounding a disk of area~$0 < a_2 < \frac{a_1}{k}$ which does not contain the orbifold point (nor the puncture), see Figure~\ref{fig:1b}. A product torus~$T(x)$ for $x \in \ell_k(a_1)$, on the other hand, projects to a closed embedded curve bounding a disk containing the orbifold point. We refer to~\S\S \ref{ssec:famsympred}-\ref{ssec:torivd} for more details.

\begin{figure}
		\begin{tikzpicture}[scale=0.85]	
			\filldraw[color=black!75, fill=black!10, very thick](0,0) circle (3);
			\fill[thick, black] (0,3)  circle[radius=2pt];
			\draw (0,3.5) node{orbifold point};
			\draw (0,-3) node[cross]{};
			\draw (0,-3.5) node{puncture};
			\filldraw[black, fill=black!20, very thick] (1,2) ..controls +(-0.8,0) and + (0.2,0.7).. (-0.5,0.5) ..controls +(-0.2,-0.7) and +(0,1).. (-1.5,-1.5) ..controls +(0,-1) and +(-0.7,-0.7).. (1.5,-1.5) ..controls +(0.7,0.7) and +(1,0).. (1,2);
			\draw[black!75, thick] (-3,0) .. controls +(0,-0.5) and +(-0.8,0).. (0,-0.65) .. controls +(0.8,0) and +(0,-0.5) .. (3,0);
			\draw[dashed ,black!75, thick] (-3,0) .. controls +(0,0.5) and +(-0.8,0).. (0,0.65) .. controls +(0.8,0) and +(0,0.5) .. (3,0);
			\node at (-1,1.2) {$C(a_2)$};
		\end{tikzpicture}
	\caption{Reduced space corresponding to a $\ell_k(a_1)$. The curve $C(a_2)$ bounds a disk of area $a_2 > 0$ and lifts to the exotic torus $\Upsilon_k(a_1,a_2) \subset \C^3$. A similar curve bounding a disk containing the orbifold point lifts to a (Hamiltonian perturbation of) a standard product torus.}
	\label{fig:1b}
\end{figure}

\begin{remark}
The segments~$\ell_k(a_1)$ intersect the edge of~$\R_{\geqslant 0}^3$ given by the first coordinate axis.
The same construction can be carried out for any face of codimension~$\geqslant 2$ in the moment polytope~$\R_{\geqslant 0}^n$ of~$\C^n$. We do not discuss the general construction here, but instead we rely on taking products to produce tori in higher dimensions, see~\S\ref{ssec:higherdim} and~\S\ref{ssec:proofthmD}.
\end{remark}

\subsection{The invariant}
\label{ssec:invariant}

The invariant we use to distinguish Lagrangian tori is the so-called \emph{displacement energy germ}. It is a strengthened symplectic invariant of Lagrangian submanifolds derived from displacement energy via \emph{versal deformations}. Versal deformations of Lagrangian submanifolds were first introduced by Chekanov~\cite{Che96}. See~\S\ref{ssec:vd} for a detailed introduction to the displacement energy germ and~\cite{Bre20, Bre22, CheSch10, CheSch16} for further applications. Here we give a brief discussion.

The \emph{displacement energy germ} of a closed Lagrangian~$L \subset (X, \omega)$ is the function
\begin{equation}
	\ce_L \colon H^1(L;\R) \rightarrow \R_{\geqslant 0} \cup \{+ \infty\},
\end{equation} 
defined as a germ at~$0 \in H^1(L;\R)$ of the displacement energy~$e(\cdot) \in \R \cup \{+ \infty\}$ viewed as a function on the space of Lagrangians. By Weinstein's neighbourhood theorem, a neighbourhood of~$L$ in the space of Lagrangians is parametrized by~$C^1$-small closed differential one-forms on~$L$. Furthermore, the graphs of two such forms are Hamiltonian isotopic in the Weinstein neighbourhood if and only if the corresponding forms are cohomologuous. The germ of Lagrangian neighbours of~$L$ up to the germ of Hamiltonian isotopies is thus in natural bijection with~$H^1(L;\R)$. Since~$e(\cdot)$ is invariant under symplectomorphisms~$\ce_L$ is well defined. The displacement energy germ is a symplectic invariant in the sense that
\begin{equation}
	\label{eq:deginv}
	\ce_{\phi(L)}  = \ce_L \circ (\phi\vert_{L})^*, \quad
	\phi \in \Symp(X,\omega).
\end{equation}
In practice, one picks a \emph{versal deformation} of~$L$, meaning a continuous family of Lagrangians $b \mapsto v_L(b)$ with~$b \in U \subset H^1(L;\R)$ which are parametrized by their Lagrangian flux, see Definition~\ref{def:vd}. Versal deformations induce the bijection discussed above and thus, to compute~$\ce_L$, one can take the germ at the origin~$0 \in H^1(L;\R)$ of the function~$b \mapsto e(v_L(b))$. This is independent of the choice of versal deformation. 

In our case,~$L$ is diffeomorphic to a torus and thus~$H^1(L;\R) \cong \R^n$ and we obtain an invariant up to transformations~$(\phi\vert_L)^* \in \GL(n;\Z)$. Furthermore, in applications, we usually determine~$b \mapsto e(v_L(b))$ on an open dense subset of a neighbourhood of the origin~$0 \in \R^n$. This yields~$\ce_L$ on an open dense subset of~$\R^n$, which is sufficient to distinguish Lagrangians. We make the following definition.

\begin{definition}
\label{def:sim}
Let~$L,L' \subset X$ be Lagrangian tori.
We say that~$\ce_L, \ce_{L'} \colon \R^n \rightarrow \R \cup \{+ \infty\}$ are \emph{equivalent} if there is~$\Phi \in \GL(n;\Z)$ such that~$\ce_{L'} = \ce_{L} \circ \Phi$ on an open dense subset of~$\R^n$. We write~$\ce_L \sim \ce_{L'}$.
\end{definition}

\subsection{Key ideas: distinguishing tori}
\label{ssec:disttori}
Let us first discuss how to distinguish tori in~$\C^n$. Let~$L = \Upsilon_k \subset \C^3$ (see Theorem~\ref{thm:A}),~$L= \Theta_\fm \subset \C^3$ (see Proposition~\ref{prop:viannalift}), or~$L = \Theta_{\fm} \times T(c_4,\ldots,c_n) \subset \C^n$ (see the proof of Theorem~\ref{thm:D} in~\S\ref{ssec:proofthmD}). To compute~$\ce_L$, we use the \emph{instability of exotic tori}, by which we mean that most tori which are~$C^1$-close to~$L$ are equivalent to product tori in~$\C^n$ and thus not exotic themselves. We note that this remarkable property is shared by all exotic tori known to the author. More precisely, for every neighbourhood~$U \subset H^1(L;\R) \cong \R^n$ of the origin, there is an open dense subset~$U' \subset U$ such that the versal deformation~$v_L(b)$ is equivalent to a product torus for all~$b \in U'$. The displacement energy of product tori is known, see Example~\ref{ex:prodtori},
\begin{equation}
	e(T(a)) = \underline{a}, \quad
	a = (a_1,\ldots,a_n) \in \R_{>0}^n,
\end{equation}
where~$\underline{a} = \min\{a_1,\ldots,a_n\}$ denotes the minimum throughout the paper. We can compute~$\ce_L$ on an open and dense subset by determining the bijection~$\sigma \colon U' \rightarrow \sigma(U') \subset \R^n$ for which
\begin{equation}
	v_L(b) \cong T(\sigma(b)), \quad
	b \in U'.
\end{equation}

Interestingly enough, in all examples we consider, the map~$\sigma$ is piecewise linear and the linear components are all in~$\GL(n;\Z)$. Note that it is sufficient to have~$\ce_L \nsim \ce_{L'}$ to conclude~$L \ncong L'$. See~\eqref{eq:upsdegerm} for the displacement energy germ of the~$\Upsilon_k$-tori and~\eqref{eq:degliftvianna} for the case of the~$\Theta_{\fm}$-tori. 

\begin{remark}
The idea of using the instability of exotic tori to compute invariants coming from versal deformations is due to Chekanov~\cite{Che96}. See~\cite[Section 3]{CheSch10} for an explicit computation for the Chekanov torus in~$S^2 \times S^2$ and see~\cite[Section 5]{Bre22} for the case of the Chekanov torus in~$\R^4$. 
\end{remark}

To distinguish embeddings of~$L$ into a geometrically bounded symplectic manifold~$(X,\omega)$ by a Darboux chart~$\psi \colon B^{2n}(R) \rightarrow X$ (see Theorems~\ref{thm:B},~\ref{thm:C} and~\ref{thm:D}) we show that
\begin{equation}
	\label{eq:degsim}
	\ce_{\psi(L)} \sim \ce_{L},
\end{equation}
provided~$L$ is small enough. We stress that~$\ce_{\psi(L)}$ is computed from displacement energy in~$X$, whereas~$\ce_{L}$ is computed from displacement energy in~$\C^n$. If~\eqref{eq:degsim} holds, then the image under Darboux embeddings of the tori~$\Upsilon_k, \Theta_{\fm}, \Theta_{\fm} \times T(c)$ are still distinguished by the displacement energy germ.

\begin{definition}
\label{def:locality}
Let~$L \subset \C^n$ be a Lagrangian. We say that an embedding~$\psi(L)$ under a Darboux chart~$\psi$ has the \emph{locality property (for the displacement energy germ)} if~\eqref{eq:degsim} holds.
\end{definition}

The main idea to prove the locality property is to combine the \emph{instability of exotic tori} with the fact that, under certain conditions, \emph{the displacement energy of product tori is preserved under Darboux embeddings}. This goes back to Chekanov--Schlenk~\cite{CheSch16}. Let us make some definitions first. Let~$J \in \cj_{\omega}$ be an~$\omega$-tame almost complex structure. Let~$\lambda_S(X,J)$ be the infimum over areas of non-constant~$J$-holomorphic spheres in~$X$. If~$\lambda_S$ is finite, then, by Gromov compactness, it is positive and realized by a $J$-holomorphic sphere. Let~$\psi \colon B(R) \rightarrow (X,\omega)$ be a Darboux chart. We will work with~$\omega$-tame almost complex structures which agree with the push-forward of the standard complex structure~$J_0$ on~$\im \psi = \psi(B(R))$. The following quantity will play a crucial role in the proofs of Theorems~\ref{thm:B} and~\ref{thm:C}, 
\begin{equation}
	\label{eq:lambdapsi}
	\lambda_S(X,\psi)
	= 
	\sup \{ \lambda_S(X,J) \, \vert \, J \in \cj_{\omega}, J\vert_{\im \psi} = \psi_* J_0 \}.
\end{equation}

\begin{definition}
\label{def:CS}
We say that a Lagrangian torus~$L \subset (X,\omega)$ has \emph{property (CS)} if there is a Darboux chart~$\psi \colon B(R) \rightarrow (X,\omega)$ with~$L = \psi(T(a))$ such that
\begin{equation}
	\label{eq:CS1}
	a_1 + \ldots + a_n + \underline{a} < R,
\end{equation}
and 
\begin{equation}
	\label{eq:CS2}
	\underline{a} < \lambda_S(X,\psi). 
\end{equation}
\end{definition}
In the aspherical case, the following proposition goes back to~\cite[Proposition 2.1]{CheSch16}. 

\begin{proposition}[Chekanov--Schlenk]
\label{prop:CS}
Let~$L = \psi(T(a)) \subset X$ be a Lagrangian torus with property {\rm (CS)}, then 
\begin{equation}
	e_X(L) = e_{\C^n}(T(a)) = \underline{a}.
\end{equation}
\end{proposition}

See~\S\ref{ssec:displacementenergy} for a proof. The following is a general criterion for when the locality property~\eqref{eq:degsim} holds.

\begin{lemma}
\label{lem:localitycrit}
Let~$L' \subset B(R)$ be a Lagrangian torus and~$\psi \colon B(R) \rightarrow (X,\omega)$ a Darboux chart. If there is a versal deformation~$v_{L'}$ defined on~$U$ such that~$v_{L'}(b) \cong T(\sigma(b))$ inside~$B(R)$ for all~$b$ in an open dense subset~$U' \subset U$ and such that all~$\psi(T(\sigma(b)))$ have property {\rm (CS)}, then the image~$L = \psi(L')$ has the locality property. 
\end{lemma}

\proof 
We identify~$H^1(L;\R)$ and~$H^1(L';\R)$ via~$\psi$. The versal deformation~$v_{L'}$ induces a natural versal deformation~$v_L = \psi \circ v_{L'}$ of~$L = \psi(L')$ and~$\ce_{L}$ is the germ of~$b \mapsto e_X(v_L(b))$. For all~$b \in U'$, we find
\begin{equation}
	e_X(v_L(b))
	= e_X(\psi(v_{L'}(b)))	
	= e_X(\psi(T(\sigma(b))))
	= e_{\C^n}(T(\sigma(b)))
	= e_{\C^n}(v_{L'}(b)).
\end{equation}
The third equality follows from Proposition~\ref{prop:CS}. Since~$\ce_{L'}$ is the germ of~$b \mapsto e_{\C^n}(v_{L'}(b))$, this proves~$\ce_L \sim \ce_{L'}$ and therefore the locality property. 
\proofend

For a fixed Darboux chart~$\psi \colon B(R) \rightarrow X$, we prove Theorems~\ref{thm:B},~\ref{thm:C} and~\ref{thm:D} by showing that the respective families of tori contain members which are small enough such that their image under~$\psi$ satisfies the hypotheses of Lemma~\ref{lem:localitycrit}. 

\begin{remark}
\label{rk:epsilon}
In Theorems~\ref{thm:B} and~\ref{thm:C}, the constant~$\varepsilon > 0$, determining the size of the tori which can be distinguished by our methods, is the minimum of a linear function of the size~$R$ of the Darboux chart and~$\lambda_S(X,\psi)$. These two terms in the minimum reflect the conditions~\eqref{eq:CS1} and~\eqref{eq:CS2} for the versal deformations to have property {\rm (CS)}. Note, that for fixed~$\psi$, the only input coming from the ambient manifold in the statements of Theorems~\ref{thm:B} and~\ref{thm:C} is the constant~$\lambda_S(X,\psi)$.
\end{remark}

\begin{remark}
\label{rk:sympaspherical}
If~$(X,\omega)$ is symplectically aspherical, we can refine the statements of Theorems~\ref{thm:B},~\ref{thm:C} and~\ref{thm:D}. In this case~$\lambda_S(X,\psi) = \infty$ and hence we can pick a sequence of ball embeddings whose sizes tend towards the Gromov width~$w = w_G(X,\omega)$ to obtain~$\varepsilon = \frac{w}{2}$ in Theorem~\ref{thm:B} and~$\varepsilon = \frac{3w}{4}$ in Theorem~\ref{thm:C}. Similarly, the statement of Theorem~\ref{thm:D} applies to all tori~$\psi(\Theta_{\fm}^n(a,a_4,\ldots,a_n))$ with~$\frac{4a}{3} + a_4 + \ldots + a_n < w$. For details, see the proofs in~\S\ref{ssec:proofthmB},~\S\ref{ssec:proofthmC} and~\S\ref{ssec:proofthmD}, respectively.
\end{remark}

\begin{remark}
\label{rk:aspherical}
Let~$X$ be aspherical. If~$L \subset \C^n$ is monotone, then so is its image~$\psi(L) \subset X$ under the Darboux chart. In that case we obtain one-parameter families of monotone exotic tori in~$X$ from Theorems~\ref{thm:B},~\ref{thm:C} and~\ref{thm:D}. 
\end{remark}

\begin{remark}
\label{rk:toricrefinement}
If~$X$ is toric, it contains a natural family of Lagrangian tori appearing as free orbits of the Hamiltonian~$T^n$-action, or equivalently, regular fibres of the moment map. These tori are called \emph{toric fibres}. We call a Lagrangian torus~$L \subset X$ \emph{exotic} if it is not equivalent to any toric fibre. In case~$X$ is monotone toric and the displacement energy of its toric fibres is, on an open and dense subset of the moment polytope, given by the integral distance to the boundary of its moment polytope, then~$X$ contains infinitely many distinct exotic tori, see Corollary~\ref{cor:toricexotic}. The condition on displacement energies is conjecturally true for all monotone toric manifolds and has been checked for all~$X$ with~$\dim X \leqslant 18$, see Example~\ref{ex:detoricfibres} or~\cite[\S 3.4]{Bre20} for more details
\end{remark}

\begin{remark}
\label{rk:cstori}
As discussed in Example~\ref{ex:cstori}, every Lagrangian torus~$L = \psi(T(a)) \subset X$ which has property {\rm (CS)} has displacement energy germ~$\ce_L \sim \ce_{T(a)}$. Having the same displacement energy germ as product tori, we may consider~{\rm (CS)}-tori to be the \emph{standard tori} in general geometrically bounded symplectic manifolds against which to compare exotic tori. Indeed, as a corollary of the proofs of Theorems~\ref{thm:B} and~\ref{thm:C}, we obtain that neither the~$\Upsilon_k$-, nor the~$\Theta_{\fm}$-tori are equivalent to a torus having property~{\rm (CS)}. Thus we recover an absolute notion of what exoticity means. 
\end{remark}

\subsection{Discussion and further questions} \label{ssec:dandq} 
In geometrically bounded symplectic manifolds of dimension at least six, the existence of an infinite collection of tori which are area-equivalent but pairwise distinct is a \emph{local} phenomenon, see Theorem~\ref{thm:D}. Relative to the current state of knowledge, this is in contrast to the situation in dimension four. As of now - meaning as long as no other tori besides product and Chekanov are found in~$\C^2$ - the existence of infinitely many exotic tori in four-dimensional symplectic manifolds is a \emph{global} phenomenon in the sense that it is highly dependent on the geometry of the ambient space. As illustrated in Example~\ref{ex:cheknonexotic}, even embedding the Chekanov torus into a manifold by a Darboux chart does not yield a torus which is exotic in a robust way. However, in dimension four there is a semi-local analogue of the constructions discussed in this paper. Instead of using Darboux charts, i.e.\ ball embeddings, one can use symplectic embeddings of more complicated domains into~$X^4$ in order to construct exotic Lagrangian tori. This is ongoing work joint with Johannes Hauber and Joel Schmitz. See~\S\ref{ssec:dimensionfour} for more details.

The second main point of this paper is that it sheds some light on the (un-)importance of monotonicity when it comes to the study of exotic tori. Most exotic tori which have been studied are monotone, see~\cite{Aur15, Bre20, Che96, CheSch10, GalMik22, Via17, Via16}, with the notable exceptions of~\cite{FOOO12} and~\cite{SheTonVia19}. However, we make the case that this circumstance is mainly an artifact of the methods used to distinguish Lagrangians. In contrast to the superpotential counting Maslov-two~$J$-holomorphic disks, the displacement energy germ does not require any adaptations in order to yield a well-defined invariant in the case of non-monotone Lagrangians. In some sense, the tori studied in Theorems~\ref{thm:B}, \ref{thm:C} and~\ref{thm:D} are the very opposite of monotone, since they bound very small Maslov two disks inside the Darboux chart as opposed to potentially very large disks coming from the ambient space. The fact that the embedded tori in Theorems~\ref{thm:B}, \ref{thm:C} and~\ref{thm:D} are chosen small relative to the size of the ball embedding contrasts with the phenomenon of exotic tori \emph{untwisting} when they are given enough space, see for example~\cite[\S 4]{CheSch10}.

An interesting direction of inquiry raised in~\S\ref{ssec:atfshigherdim} is that of interpreting the~$\Upsilon_k$-tori as fibres of Lagrangian torus fibrations on~$\C^3$ arising from mutations of its toric structure, just as the Vianna tori in~$\CP^2$ appear as fibres of almost toric mutations of the toric structure on~$\CP^2$. This seems to lead to an analogue of almost toric fibrations in higher dimensions which would provide a framework to carry out higher mutations in arbitrary dimensions to produce more exotic tori. Taking fibres of such mutated Lagrangian torus fibrations should provide a general framework in which both the~$\Upsilon_k$- and the~$\Theta_{\fm}$-tori appear as special cases. 

\subsection{Structure of the paper} Section~\ref{sec:background} provides background on symplectic reduction, integral geometry, displacement energy and the displacement energy germ which is used in later sections. Section~\ref{sec:atfs} gives an overview of \emph{almost toric fibrations}, mutations of ATFs on~$\CP^2$ and a computation of the displacement energy germ of certain almost toric fibres. The approach we follow in~\S\ref{ssec:atfibres} was suggested to us by Felix Schlenk based on discussions with Yuri Chekanov. The displacement energy germ of almost toric fibres in~$\CP^2$ is used to distiguish the tori~$\Theta_{\fm}(a) \subset \C^3$ and their embeddings, see Proposition~\ref{prop:viannalift} and Theorems~\ref{thm:C},~\ref{thm:D}. The heart of the paper can be found in Section~\ref{sec:constrproofs}, which contains the construction of the~$\Upsilon_k$-tori, as well as all proofs. In Section~\ref{sec:qandp}, we discuss open questions, we give an outlook on future work in dimension four and speculate on mutations of Lagrangian torus fibrations in dimensions greater than or equal to six, which may lead to a theory of almost toric fibrations in these dimensions.

\subsection*{Acknowledgements}
It is our pleasure to thank Jonny Evans, Jeff Hicks, Yael Karshon, Leonid Polterovich and Joel Schmitz for interesting discussions, encouragement and insightful comments on an earlier draft. We thank Denis Auroux for a stimulating exchange about the feasability of using the count of Malsov index two $J$-holomorphic disks to obtain similar results. We thank Renato Vianna for sharing insights into his discussions joint with Jeff Hicks and Dmitry Tonkonog, see Remark~\ref{rk:vianna}. We thank Felix Schlenk for detailed comments on an earlier draft and for generously sharing his ideas joint with Yuri Chekanov, in particular those surrounding~\S\ref{ssec:atfibres}. We are endebted to Paul Biran for a crucial remark pointing out a mistake after the author's talk at ETH in October 2022. We thank the anonymous referee for many useful remarks and in particular, for pointing out the relevance of \cite{HinZha23} for Question \ref{q:1}, see Remark \ref{rk:stablyexotic}. We acknowledge the support of the following grants: Israel Science Foundation grant 1102/20, ERC Starting Grant 757585 and Swiss National Science Foundation Ambizione Grant PZ00P2-223460.

\section{Background}
\label{sec:background}

\subsection{Symplectic reduction}
\label{ssec:sympred}
In this subsection we introduce some notions surrounding symplectic reduction which will be used in subsequent sections. Although most properties hold in the more general context of Hamiltonian actions of compact Lie groups, we restrict our attention to Hamiltonian torus actions. We refer to~\cite{Can01, GuiSte90, McDSal17, Mei00} for more details on this classical topic. Let~$(X,\omega)$ be a symplectic manifold equipped with a Hamiltonian~$T^k$-action, generated by the moment map~$\mu \colon X \rightarrow \R^k$. If the restriction of the~$T^k$-action to a level set~$\mu^{-1}(c)$ is free, then the corresponding level set is smooth and its quotient~$X_c = \mu^{-1}(c)/T^k$ is a smooth manifold equipped with a naturally induced symplectic form~$\omega_c$ satisfying~\eqref{eq:reductionomega}. This procedure is called \emph{symplectic reduction} and the symplectic manifold~$(X_c,\omega_c)$ is called \emph{symplectic quotient}. We will use a slightly generalized version of symplectic reduction, where we allow for group actions having finite stabilizers. In that case, the reduced space is a symplectic orbifold, see for example~\cite{LerTol97} for details. Symplectic reduction is best summarized by a diagram, 
\begin{equation}
\label{eq:reductiondiag}
\begin{tikzcd}
	\mu^{-1}(c) \arrow[r, hook] \arrow[d, "\pi"] &
	(X,\omega) \\´
	(X_c, \omega_c),
\end{tikzcd}
\end{equation}
with
\begin{equation}
	\label{eq:reductionomega}
	\pi^* \omega_c = \omega\vert_{\mu^{-1}(c)}.
\end{equation}
A crucial feature of symplectic reduction we will use repeatedly in this paper is the fact that one can lift Lagrangian submanifolds and Hamiltonian isotopies in a natural way from the reduced space~$X_c$ to~$X$, see for example~\cite{ AbrMac13, Bre20, Bre23}.

\begin{proposition}
\label{prop:redlifting}
Let~$L \subset X_c$ be a Lagrangian, then~$\pi^{-1}(L) \subset \mu^{-1}(c) \subset X$ is a~$T^k$-invariant Lagrangian. Let~${\varphi_t}$ be a Hamiltonian isotopy of~$X_c$, then there is a~$T^k$-equivariant Hamiltonian isotopy~$\{\phi_t\}$ satisfying~$\pi \circ \phi_t = \varphi_t \circ \pi$. 
\end{proposition}

Similarly, invariant Lagrangians and equivariant Hamiltonian isotopies in~$X$ project to the symplectic quotient~$X_c$. For Hamiltonian group actions we have the following.  

\begin{proposition}
\label{prop:redgroupaction}
Let~$\nu \colon X \rightarrow \R^l$ be another Hamiltonian~$T^l$-action which commutes with the~$T^k$-action, i.e.\ $\{\mu_i,\nu_j\} = 0 $ for all~$i,j$. Then there is an induced Hamiltonian~$T^l$-action on the symplectic quotient~$X_c$ generated by~$\nu_c \colon X_c \rightarrow \R^l$ with~$\pi^* \nu_c = \nu\vert_{\mu^{-1}(c)}$ for which~$\pi$ is~$T^l$-equivariant.
\end{proposition}

A special case of this is \emph{toric reduction}. A toric manifold is a symplectic~$2n$-manifold~$(X, \omega)$ equipped with an effective Hamiltonian~$T^n$-aciton. The image of the moment map~$\Delta = \mu(X) \subset \R^n$ is a polytope with special properties, a so-called \emph{Delzant polytope}, whose integral affine geometry classifies~$(X,\omega,\mu)$ up to equivariant symplectomorphisms by the Delzant Theorem~\cite{Del88}. We refer to~\cite{Eva23, Gui94} for more details. Performing symplectic reduction with respect to a subtorus~$K \subset T^n$ of the toric action, we are in the context of Proposition~\ref{prop:redgroupaction}, since the actions of~$K$ and~$T^n$ commute. The quotient~$X_c$ carries an effective~$T^n/K$-action and is thus a toric manifold in its own right. On the level of moment polytopes, toric reduction is encoded by an integral affine embedding of the Delzant polytope of~$X_c$ into~$\Delta$. The Delzant construction of toric manifolds~\cite[Chapter 1]{Gui94} and McDuff's probes~\cite{McD11} are special cases of toric reduction. The symplectic reductions in~\S\ref{ssec:famsympred} are of this type, too, and can be seen as performing toric reduction on a family of segments in the moment polytope~$\R^3_{\geqslant 0}$ of~$\C^3$. We refer to~\cite[\S 2.2]{Bre23} for more details. 

Regular fibres~$\mu^{-1}(x)$ of toric moment maps are Lagrangian tori, see also Example~\ref{ex:degtoricfibres}. As a special case of Proposition~\ref{prop:redlifting}, recall from~\cite[Proposition 2.8]{Bre23} that toric fibres lift to toric fibres under the toric reduction. 

\begin{proposition}
\label{prop:redtoricfibres}
Let~$(X_c,\omega_c,\mu_c)$ be the toric reduction of a toric~$(X,\omega,\mu)$ with~$\iota \colon \Delta_c \hookrightarrow \Delta$ the corresponding integral affine embedding of Delzant polytopes. Then
\begin{equation}
	\mu^{-1}(\iota(x)) = \pi^{-1}(\mu_c^{-1}(x)), \quad
	x \in \Int \Delta_c. 
\end{equation}
\end{proposition}

\subsection{Some integral geometry}

Integral geometry is the geometry associated to transformations preserving a lattice of maximal rank in a finite-dimensional vector space. We restrict our attention to the lattice~$\Z^n \subset \R^n$ which has symmetry group~$\GL(n;\Z)$. Our interest in integral geometry comes from the fact that, as discussed in~\S\ref{ssec:invariant} and~\S\ref{ssec:vd}, the displacement energy germ of Lagrangian tori is an invariant up to~$\GL(n;\Z)$. For related reasons, the regular locus of the base of any Lagrangian torus fibration has an \emph{integral affine structure}. This means that the Lagrangian torus fibration induces a manifold atlas with transition charts in~$\GL(n;\Z) \ltimes \R^n$, which essentially comes from action coordinates and the Arnold--Liouville theorem, see~\cite[\S 2.3]{Eva23}. The integral affine structure of Lagrangian torus fibrations is intimately related to the versal deformations of its regular fibres, see Proposition~\ref{prop:momentmapsvd}.

A vector~$v \in \Z^n$ is called \emph{primitive} if its entries are coprime. Let~$w \in \R^n$ be a vector with a rational direction. Then there is a unique~$\alpha > 0$ such that~$w = \alpha v$ for a primitive~$v \in \Z^n$ and~$\alpha$ is called \emph{integral affine length of~$w$} and denoted by~$\ell_{\rm IA}(w) = \alpha$. Note that this corresponds to the Euclidean length of the image of~$w$ under a map in~$\GL(n;\Z)$ mapping~$w$ into a coordinate axis of~$\R^n$. Let~$P \subset \R^n$ be a rational affine hyperplane. Then the \emph{integral affine distance}~$d_{\IA}(x,P)$ of a point~$x$ to~$P$ is defined as the Euclidean distance of~$\phi(x)$ to~$\{x_n = 0\}$, where~$\phi \in \GL(n,\Z) \ltimes \R^n$ is an integral affine transformation with~$\phi(P) = \{x_n = 0\}$. This is independent of the choice of~$\phi$. Write~$P = \{y \in \R^n \, \vert \, \langle y , v \rangle + \lambda = 0 \}$ for a primitive~$v \in \Z^n$ and~$\lambda \in \R$, then~$d_{\rm IA}(x,P) = \vert \langle x, v \rangle + \lambda \vert$. Now let
\begin{equation}
	\label{eq:momentpolytope}
	\Delta = \{ x \in \R^n \;\vert\; \ell_i(x) \geqslant 0 \}, \quad
	\ell_i(x) = \langle x , v_i \rangle + \lambda_i
\end{equation}
be a convex polytopes defined by the intersection of half-spaces delimited by rational affine hyperplanes~$P_i = \{\ell_i(y) = 0\}$ for~$i \in \{1, \ldots, N\}$. Again, the~$v_i \in \Z^n$ are chosen to be primitive. They can be thought of as the inwards pointing normal vectors to the facets of~$\Delta$. Then the \emph{integral affine distance}~$d_{\rm IA}(x,\pp \Delta)$ of~$x \in \Delta$ to the boundary~$\pp \Delta$ is defined as 
\begin{equation}
	\label{eq:iadistanceboundary}
	d_{\rm IA}(x, \pp \Delta) 
	= 
	\min_{1 \leqslant i \leqslant N}\{d_{\rm IA}(x,P_i)\}
	=
	\min\{ \ell_1(x), \ldots, \ell_N(x)\}.
\end{equation}

\begin{definition}
\label{def:intindex}
Let~$v_1,\ldots,v_k \in \Z^n$ for~$k \leqslant n$. The \emph{integral index}~$D_{\rm IA} = D_{\rm IA}(v_1,\ldots,v_k)$ is defined as the greatest common divisor of
\begin{equation}
	\left\{ \det(v_1 \vert \ldots \vert v_k \vert w_{k+1} \vert \ldots \vert w_n) \in \Z \,\vert\,
	 w_{k+1}, \ldots, w_n \in \Z^n \right\}.
\end{equation}
\end{definition} 

This is an integral invariant of the unordered~$k$-tuple~$v_1,\ldots,v_k$. Its geometric interpretation is that it is equal to the index of the sublattice~$\Span_{\Z}\{v_1,\ldots,v_k\}$ of the lattice~$\Z^n \cap \Span_{\R}\{v_1,\ldots,v_k\}$. For computations, the following characterization is the most convenient, see~\cite[Chapter 18]{Kar22} for a proof and more details.

\begin{proposition}
\label{prop:intindex}
The integral index~$D_{\rm IA}(v_1,\ldots,v_k) \in \Z$ is equal to the greatest common divisor of all~$k$-minors (meaning determinants of~$k \times k$ submatrices) of the matrix~$(v_1\vert \ldots \vert v_k)$.
\end{proposition}

\subsection{Displacement energy}
\label{ssec:displacementenergy}
Let~$A$ be a compact subset of a symplectic manifold~$(X,\omega)$. Recall that the \emph{displacement energy of~$A$} is defined as
\begin{equation}
\label{eq:displacementenergy}
	e(A) 
	= 
	\inf \{ \Vert H \Vert \,\vert\, \phi^H_1(A) \cap A = \varnothing \},
\end{equation}
where~$\phi^H_1$ denotes the time-one flow of a compactly supported, time-dependent Hamiltonian~$H$. We set~$e(A) = \infty$ if the infimum in~\eqref{eq:displacementenergy} is taken over the empty set and by~$\Vert \cdot \Vert$ we denote the Hofer norm,
\begin{equation}
	\Vert H \Vert 
	= 
	\int_0^1 \left( \max_{x \in X} H_t(x) - \min_{x \in X} H_t(x) \right) dt.
\end{equation}
We refer to the classic~\cite{Pol01} for details on Hofer geometry. The quantity~$e(A)$ is a symplectic invariant of the set~$A$, in the sense that if there is a symplectomorphism of~$(X,\omega)$ mapping~$A$ to another subset~$B \subset X$, then the displacement energies of~$A$ and~$B$ agree. In this article, we are exclusively interested in the displacement energy of Lagrangian tori. Chekanov~\cite{Che98} (see also~\cite{Oh97}) gave a lower bound for displacement energy of compact Lagrangian submanifolds in geometrically bounded symplectic manifolds in terms of~$J$-holomorphic spheres in~$X$ and~$J$-holomorphic disks with boundary on the Lagrangian in question. Denote by~$\cj = \cj(X,\omega)$ the set of~$\omega$-tame almost complex structures. Let~$\lambda_D(X,L,J)$ denote the infimum of symplectic areas realized by non-trivial~$J$-disks with boundary on~$L$. If such disks exist, then the infimum is positive and is attained by at least one disk by Gromov compactness. If not, then set~$\lambda_D(X,L,J) = \infty$. Similarly, denote by~$\lambda_S(X,J)$ the infimum of areas realized by non-trivial~$J$-spheres.

\begin{theorem}[Chekanov]
\label{thm:Chekanov}
Let~$L \subset X$ be a compact Lagrangian in a geometrically bounded symplectic manifold and~$J$ an~$\omega$-tame almost complex structure. Then
\begin{equation}
	\min\{ \lambda_D(X,L,J) , \lambda_S(X,J) \} \leqslant e(L).
\end{equation}
\end{theorem}

\begin{example}[Product tori in~$\C^n$]
\label{ex:prodtori}
Let~$\C$ be equipped with the standard symplectic form and let~$S^1(x) \subset \C$ be the circle enclosing symplectic area~$x > 0$. Then~$e(S^1(x)) = x$. The upper bound comes from an explicit construction and the lower bound comes from the obvious holomorphic disk with boundary on~$S^1(x)$ together with Theorem~\ref{thm:Chekanov}. Taking products of this example, we obtain the \emph{product tori} in~$\C^n$, which we will denote by
\begin{equation}
	\label{eq:prodtori}
	T(x) 
	= 
	S^1(x_1) \times \ldots \times S^1(x_n) \subset \C^n, \quad
	x = (x_1,\ldots,x_n) \in \R^n_{>0}.
\end{equation}
Product tori have displacement energy
\begin{equation}
	\label{eq:deprodtori}
	e(T(x))
	=
	\min\{x_1,\ldots,x_n\}.
\end{equation}
This was proved by Sikorav~\cite{Sik90} without using~$J$-holomorphic curves. It can also be proved using Theorem~\ref{thm:Chekanov} and the standard holomorphic structure on~$\C^n$.  

\end{example}

\begin{example}[Toric fibres]
\label{ex:detoricfibres}
Let~$(X,\omega)$ be a toric manifold with moment map~$\mu$ and moment polytope~$\Delta = \{x \in \R^n \, \vert \, \ell_i(x) \geqslant 0\}$, with the same conventions as in~\eqref{eq:momentpolytope}. Every fibre~$T(x) = \mu^{-1}(x)$ over a point~$x \in \Int \Delta$ in the interior of~$\Delta$ is a Lagrangian torus, called \emph{toric fibre}. In all examples known to the author,
\begin{equation}
	\label{eq:detoricfibres}
	e(T(x)) 
	= 
	d_{\rm IA}(x , \pp \Delta)
	=
	\min_{1 \leqslant i \leqslant N}\{\ell_i(x)\},
\end{equation}
for all~$x$ in an open dense subset of~$\Int \Delta$, where~$d_{\rm IA}(\cdot, \pp \Delta)$ denotes the integral affine distance to the boundary of~$\Delta$, see~\eqref{eq:iadistanceboundary}. See Figure~\ref{fig:2} for an illustration of the level sets of the function $x \mapsto d_{\rm IA}(x, \pp \Delta)$. Note that Example~\ref{ex:prodtori} is a special case of this, since product tori are the fibres of the standard toric moment map on~$\C^n$. There are toric fibres where~\eqref{eq:detoricfibres} does not hold; e.g.\ Entov--Polterovich~\cite{EntPol06} proved that every compact toric manifold contains at least one non-displaceable fibre, see also~\cite{FOOO10, FOOO11, FOOO12b}. Note that for example in~$\CP^2$, the non-displaceable fibre is the only point where~\eqref{eq:detoricfibres} fails. The inequality~$e(T(x)) \geqslant d_{\rm IA}(x, \pp \Delta)$ was shown in~\cite[Proposition 3.2]{Bre20} for \emph{all} toric fibres in compact toric manifolds. The upper bound on displacement energy relies on probes~\cite{McD11} and therefore on delicate properties of the moment polytope. If~$(X^{2n},\omega)$ is monotone toric, the upper bound follows from a conjecture (related to the Ewald conjecture, see~\cite[\S 3.1]{McD11}) about lattice polytopes, which has been checked in dimensions~$2n \leqslant 18$, see~\cite[\S 3.4]{Bre20} for more details.
\end{example}

\begin{figure}
		\begin{tikzpicture}[scale=0.55]	
			\fill[black!10] (5,5)--(5,-3)--(3,-5)--(-5,-5)--(-5,-1)--(-3,3)--(-1,5)--(5,5);
			\draw[thick,black] (5,5)--(5,-3)--(3,-5)--(-5,-5)--(-5,-1)--(-3,3)--(-1,5)--(5,5);
			\draw[thick,black!50] (4.5,4.5)--(4.5,-3)--(3,-4.5)--(-4.5,-4.5)--(-4.5,-0.5)--(-3,2.5)--(-1,4.5)--(4.5,4.5);
			\draw[thick,black!50] (4,4)--(4,-3)--(3,-4)--(-4,-4)--(-4,0)--(-3,2)--(-1,4)--(4,4);
			\draw[thick,black!50] (3.5,3.5)--(3.5,-3)--(3,-3.5)--(-3.5,-3.5)--(-3.5,0.5)--(-3,1.5)--(-1,3.5)--(3.5,3.5);
			\draw[thick,black!50] (3,3)--(3,-3)--(-3,-3)--(-3,1)--(-1,3)--(3,3);
			\draw[thick,black!50] (2.5,2.5)--(2.5,-2.5)--(-2.5,-2.5)--(-2.5,1)--(-1,2.5)--(2.5,2.5);
			\draw[thick,black!50] (2,2)--(2,-2)--(-2,-2)--(-2,1)--(-1,2)--(2,2);
			\draw[thick,black!50] (1.5,1.5)--(1.5,-1.5)--(-1.5,-1.5)--(-1.5,1)--(-1,1.5)--(1.5,1.5);
			\draw[thick,black!50] (1,1)--(1,-1)--(-1,-1)--(-1,1)--(1,1);
			\draw[thick,black!50] (0.5,0.5)--(0.5,-0.5)--(-0.5,-0.5)--(-0.5,0.5)--(0.5,0.5);
			\draw[dotted,thick,black!80] (0,0)--(5,5);
			\draw[dotted,thick,black!80] (0,0)--(3,-3);
			\draw[dotted,thick,black!80] (3,-3)--(5,-3);
			\draw[dotted,thick,black!80] (3,-3)--(3,-5);
			\draw[dotted,thick,black!80] (0,0)--(-5,-5);
			\draw[dotted,thick,black!80] (0,0)--(-1,1);
			\draw[dotted,thick,black!80] (-1,1)--(-3,1);
			\draw[dotted,thick,black!80] (-1,1)--(-1,5);
			\draw[dotted,thick,black!80] (-3,1)--(-3,3);
			\draw[dotted,thick,black!80] (-3,1)--(-5,-1);
			\fill[thick, black!80] (0,0)  circle[radius=1.5pt];
		\end{tikzpicture}
	\caption{A (non-monotone) Delzant polytope and level sets of the integral affine distance to the boundary. Chambers $\{ d_{\rm IA}(x , \pp \Delta) = \ell_i(x)\}$ are delimited by dotted lines. On an open dense subset, the integral affine distance corresponds to the displacement energy.}
	\label{fig:2}
\end{figure}

Let us prove Proposition~\ref{prop:CS} along the lines of~\cite{CheSch16}. \smallskip

\proofof{Proposition~\ref{prop:CS}}
Let~$L = \psi(T(a)) \subset X$ be a Lagrangian torus which has property {\rm (CS)}. For the inequality~$e_X(L) \leqslant \underline{a}$, see~\cite[Proposition 2.1]{CheSch16} or consider the standard moment map~$\mu \colon \C^n \rightarrow \R_{>0}^n$ to view~$B(R)$ as an open toric manifold and use McDuff's method of probes~\cite{McD11}. The condition~\eqref{eq:CS1} precisely guarantees that there is a sufficiently long probe contained in the moment map image of~$B(R)$ to displace~$T(a)$. For the inequality~$e_X(L) \geqslant \underline{a}$, let~$\{J_k\}_{k \in \N}$ be a sequence of~$\omega$-tame almost complex structures satisfying~$J_k \vert_{\im \psi} = \psi_* J_0$ such that~$\lambda_S(X,J_k) \rightarrow \lambda_S(X,\psi)$ for~$k \rightarrow + \infty$. By Theorem~\ref{thm:Chekanov}, we have 
\begin{equation}
	\label{eq:ChekanovJk}
	\min \{ \lambda_D(X,L,J_k), \lambda_S(X,J_k) \} \leqslant e(L), \quad
	k \in \N.
\end{equation}
It follows from the proof of~\cite[Proposition 2.1]{CheSch16} that~$\lambda_D(X,L,J_k) = \underline{a}$. The essential ingredient is the fact that the minimal~$\psi_*J_0$-holomorphic disk inside~$\im \psi$ has area~$\underline{a}$ and that all disks leaving~$\im \psi$ have area~$> \underline{a}$, which follows from the area estimate~\cite[Appendix A]{CheSch16} together with~\eqref{eq:CS1}. Since~$\lambda_S(X,J_k) \rightarrow \lambda_S(X,\psi)$ and~$\underline{a} < \lambda_S(X,\psi)$ by~\eqref{eq:CS2}, there is~$k_0 \in \N$ such that the left hand side of~\eqref{eq:ChekanovJk} is equal to~$\underline{a}$. 
\proofend

\subsection{Versal deformations and the displacement energy germ}

\label{ssec:vd}

We give a brief overview of the method of versal deformations here. When applied to displacement energy, it yields the \emph{displacement energy germ}, which is the invariant used in this paper to distinguish Lagrangian submanifolds. See also~\cite{Che96, CheSch16} and~\cite{Bre20, Bre22, CheSch10} for further applications.

Let~$L \subset X$ be a compact Lagrangian and let~$\cl = \cl_L$ be the space of Lagrangians which are Lagrangian isotopic to~$L$, equipped with the~$C^1$-topology, see for example~\cite{Ono07} for details. In this topology, a small enough neighbourhood~$\cu$ of~$L \in \cl$ is in bijection with a neighbourhood of the zero $1$-form in the set of closed differential one-forms on~$L$. This follows from Weinstein's neighbourhood theorem applied to~$L$ together with the fact that a graph~$\Gamma_{\alpha} \subset (T^*L,d\lambda)$ of a differential $1$-form~$\alpha \in \Omega^1(L)$ is Lagrangian if and only if~$\alpha$ is closed. Furthermore, the graphs~$\Gamma_{\alpha},\Gamma_{\beta}$ are Hamiltonian isotopic in~$T^*L$ if and only if the corresponding forms are cohomologous,~$[\alpha] = [\beta] \in H_{\rm dR}^1(L;\R)$. In other words, there is a continuous bijection between~$C^1$-small Lagrangian perturbations of~$L$ up to locally supported Hamiltonian isotopies and small elements in~$H^1(L;\R)$. In terms of a Weinstein chart~$\varphi \colon T^*L \supset V \rightarrow X$, this bijection can be realized by taking the class up to Hamiltonian isotopies supported in~$\im \varphi$ of~$v_L^{\varphi} (b) = \varphi(\Gamma_{\beta(b)})$. Here~$b \mapsto \beta(b)$ is a continuous family of closed $1$-forms satisfying~$[\beta(b)] = b$. 

\begin{example}
Let~$L \subset X$ be a Lagrangian torus and~$\varphi \colon T^*T^n \supset V \rightarrow X$ be a Weinstein chart. Identifying~$T^n = (\R/\Z)^n$, we obtain induced identifications~$T^*T^n \cong (\R / \Z)^n \times \R^n = \{(\theta_1,\ldots,\theta_n,p_1,\ldots,p_n)\}$ and~$H^1(T^n;\R) \cong \R^n$. Set~$\beta(b) = \sum_i b_i d\theta_i$ to obtain
\begin{equation}
	\label{eq:toriweinsteinvd}
	v^{\varphi}_L(b) = \varphi(T^n \times \{b\}).
\end{equation}
\end{example}

The above discussion is closely related to the Lagrangian flux of Lagrangian isotopies contained in the image of~$\varphi$. We recall the definition of Lagrangian flux for the reader's convenience, see also~\cite{Ono07}. Let~$\{i_t\}_{t \in [0,1]}$ be a Lagrangian isotopy, meaning that there is a smooth map~$i \colon [0,1] \times L \rightarrow X$ such that~$i_t = i(t,\cdot)$ is a Lagrangian embedding for every~$t \in [0,1]$. The embeddings being Lagrangian is equivalent to the forms~$i_t^*\iota(X_t)\omega \in \Omega^1(L)$ being closed for all~$t \in [0,1]$, where~$X_t \in \Gamma(i_t^*TX)$ is the vector field generating~$i_t$. The \emph{Lagrangian flux of~$\{i_t\}$} is defined as the average cohomology class of this family of forms,
\begin{equation}
	\label{eq:lagflux}
	\Flux \{i_t\} = \int_0^1 [i_t^*\iota(X_t)\omega] \in H^1(L;\R). 
\end{equation}
As in the case of the flux of symplectic isotopies, see~\cite[Chapter 14]{Pol01} or~\cite[\S 10.2]{McDSal17}, the flux of Lagrangian isotopies can be geometrically interpreted as a measure of the symplectic area of cylinders swept out by one-cycles. Let~$c$ be a one-cycle on~$L$ with homology class~$\xi = [c] \in H_1(L)$. We obtain a two-chain~$C_{\xi} = \cup_{t \in [0,1]} i_t(c)$ in~$X$, defining a relative homology class~$[C_{\xi}] \in H_2(X,\iota_0(L) \cup \iota_1(L))$. The assignment~$\xi \mapsto [C_{\xi}]$ yields a well-defined, linear map~$H_1(L) \rightarrow H_2(X,\iota_0(L) \cup \iota_1(L))$. Since~$\iota_0,\iota_1$ are Lagrangian embeddings, the symplectic form yields a well-defined functional~$[\omega] \colon H_2(X,\iota_0(L) \cup \iota_1(L)) \rightarrow \R$. From this point of view, Lagrangian flux can be viewed as
\begin{equation}
	\label{eq:lagfluxalt}
	\langle \Flux \{i_t\} , \xi \rangle
	= 
	\langle [\omega] , [C_{\xi}] \rangle, \quad
	\xi \in H_1(L).
\end{equation}
In this paper, we work with unparametrized Lagrangian isotopies, i.e.\ families~$\{L_t\}_{t \in [0,1]}$ of Lagrangian submanifolds. One can always choose a parametrization~$i_t(L) = L_t$ and prove that the above discussion is independent of this choice up to picking an identification~$L = L_0$. Furthermore, we sometimes consider \emph{continuous} families of Lagrangians, i.e.\ continuous maps~$[0,1] \rightarrow \cl$. By its homological interpretation~\eqref{eq:lagfluxalt}, the Lagrangian flux is well-defined. Furthermore since~$\omega$ is closed, Stokes' theorem and~\eqref{eq:lagfluxalt} show that the Lagrangian flux is invariant under homotopies of Lagrangian isotopies with fixed endpoints. 

\begin{remark}
A Lagrangian isotopy~$\{i_t\}$ is induced by an ambient Hamiltonian isotopy if and only if the forms~$i_t^*\iota(X_t)\omega$ are exact, i.e.\ the Lagrangian flux~$\Flux \{i_t\}_{t \in [0,t_0]}$ is zero for all~$t_0 \in [0,1]$. 
\end{remark}

\begin{remark}
\label{rk:fluxareaclasses}
The Lagrangian flux measures the difference between the area classes~$\sigma_{L_0}, \sigma_{L_1}$ (see~\eqref{eq:softinvvv}) of endpoints of a Lagrangian isotopy~$\{L_t\}$, in the sense that
\begin{equation}
	\langle \Phi^* \sigma_{L_1} - \sigma_{L_0}, D \rangle
	= 
	\langle \Flux\{L_t\} , \pp D \rangle, \quad
	D \in H_2(X,L_0),
\end{equation}
where~$\Phi \colon H_2(X,L_0) \rightarrow H_2(X,L_1)$ is the isomorphism induced by~$\{L_t\}$ and~$\pp \colon H_2(X,L_0) \rightarrow H_1(L_0)$ is the boundary map of the long exact sequence of the pair~$(X,L_0)$. In particular, a Lagrangian isotopy with vanishing flux and endpoints on non-equivalent Lagrangians yields an example answering Question~\ref{q:main} positively. The existence of such Lagrangian isotopies (see for example Remark~\ref{rk:areaequivalent}) is in contrast with the case of symplectic flux, where every isotopy with vanishing symplectic flux is homotopic to a Hamiltonian isotopy, see~\cite[Theorem 10.2.5]{McDSal17}.
\end{remark}

Let us return to the map~$v_L^{\varphi}$ defined above. Since~$\varphi^*\omega = d \lambda$ on~$\im \varphi$, Stokes' theorem shows that the flux of a Lagrangian isotopy~$\{L_t\}$ contained in~$\im \varphi$ depends only on its endpoints. Furthermore, if~$L_1 = \varphi(\Gamma_{\beta})$, then~$\Flux \{L_t\} = [\beta] \in H^1(L;\R)$. This means that the Lagrangian flux induces an inverse map to the bijection~$v_L^{\varphi}$. The map~$v_{L}^{\varphi}$ is an example of a versal deformation. However, we sometimes work with versal deformations without considering an underlying Weinstein chart~$\varphi$. In general, a versal deformation is defined as a right-inverse to the Lagrangian flux map.

\begin{definition}
\label{def:vd}
A continuous map~$v_L \colon H^1(L;\R) \supset U \rightarrow \cl_L$ defined on a simply-connected neighbourhood~$U$ of the origin is called a \emph{versal deformation} of~$L$ if~$v_L(0) =L$ and if the Lagrangians in its image are parametrized by their flux. By this we mean that
\begin{equation}
	\Flux \{ v_L(\gamma(t)) \}_{t \in [0,1]}
	= b,
\end{equation} 
for any~$b \in U$ and any continuous curve~$\gamma \colon [0,1] \rightarrow U$ with~$\gamma(0) = 0$ and~$\gamma(1) = b$.
\end{definition}

Although versal deformations are not unique, they are locally unique up to Hamiltonian isotopy. More precisely, let~$v_{L} \colon U \rightarrow \cl_L$ and~$v_L' \colon U' \rightarrow \cl_L$ be two versal deformations of~$L$. Then there is a neighbourhood of the origin~$V \subset U \cap U'$ such that
\begin{equation}
	\label{eq:vduniquene}
	v_L(b) \cong v_L'(b), \quad
	b \in V.
\end{equation}
This follows from the fact that, given a Weinstein chart~$\varphi$ of~$L$, there is a small enough neighbourhood~$V$ of the origin~$0 \in H^1(L;\R)$ such that every Lagrangian in the image of~$v_L\vert_V$ and~$v_L'\vert_V$ is the image of a graph of a closed $1$-form under~$\varphi$, and closed $1$-forms are classified up to Hamiltonian isotopy by their cohomology class. Furthermore, let~$\phi \in \Symp(X,\omega)$ and~$v_{\phi(L)}$ be a versal deformation of~$\phi(L)$, then there is a neighbourhood of the origin~$V' \subset H^1(\phi(L);\R)$ such that
\begin{equation}
	\label{eq:vdinvarian}
	v_{\phi(L)}(b) \cong (\phi \circ v_{L} \circ \phi\vert_L^*)(b), \quad
	b \in V'.
\end{equation}

Chekanov's idea is to evaluate a symplectic invariant of Lagrangian submanifolds on the image of a versal deformation and then take its germ to extract an invariant. In our case, we take displacement energy as symplectic invariant and consider the function 
\begin{equation}
	\label{eq:defct}
	H^1(L;\R) \supset U \rightarrow \R \cup \{+ \infty\}, \quad
	b \mapsto e(v_L(b)),
\end{equation}
where~$v_L$ is a versal deformation of~$L$.

\begin{definition}
Let~$L \subset (X,\omega)$ be a compact Lagrangian.
The germ~$\ce_L \colon H^1(L;\R) \rightarrow \R \cup \{+ \infty\}$ at~$0 \in H^1(L;\R)$ of the function~\eqref{eq:defct} is called \emph{displacement energy germ of~$L$}. 
\end{definition}

By~\eqref{eq:vduniquene} and the invariance of~$e(\cdot)$ under symplectomorphisms, the displacement energy germ is independent of the choice of versal deformation. The symplectic invariance property~\eqref{eq:deginv} of the displacement energy germ similarly follows from~\eqref{eq:vdinvarian}. We are interested in versal deformations of Lagrangian tori~$L \approx T^n$. It is useful to note that if~$L$ is an orbit of a Hamiltonian~$T^n$-actions, then the moment map generating the Hamiltonian~$T^n$-action yields a natural versal deformation of~$L$. Indeed, let~$\varphi \colon T^n \times U \rightarrow X$ be a~$T^n$-equivariant Weinstein chart of~$L$ and consider the versal deformation~$v_L^{\varphi} \colon U \rightarrow \cl_L$ defined by~$v_L^{\varphi}(b) = \varphi(T^n \times \{b\})$. By~$T^n$-equivariance, this proves the following. 

\begin{proposition}
\label{prop:momentmapsvd}
Let~$L \subset X$ be a Lagrangian torus and~$\mu \colon U \rightarrow \R^n$ the moment map of a (locally defined) Hamiltonian~$T^n$-action having~$L = \mu^{-1}(a)$ as orbit. Then~$v_{L}(b) = \mu^{-1}(a + b)$ is a versal deformation. 
\end{proposition}

In particular, a versal deformation of a regular fibre of a Lagrangian torus fibration (or integrable system) can be computed by computing local action-coordinates on the basis. Recall that this is equivalent to computing the flux map, see~\cite[Lemma 2.15]{Eva23}, which connects back to the appearance of Lagrangian flux in the discussion of versal deformations.

\begin{example}[Product tori in~$\C^n$]
\label{ex:degprodtori}
Let~$T(a) \subset \C^n$ be a product torus. It is easy to see that~$v_{T(a)}(b) = T(a+b)$ is a versal deformation, e.g.\ as a special case of Proposition~\ref{prop:momentmapsvd}, since the standard moment map~$\mu$ on~$\C^n$ satisfies~$\mu^{-1}(a) = T(a)$. By Example~\ref{ex:prodtori}, its displacement energy germ is the germ of~$b \mapsto \min_i\{a_i + b_i\}$ and is thus given by
\begin{equation}
	\ce_{T(a)}(b)
	=
	\underline{a} + \min_{i \in I(a)}\{b_i\},
\end{equation}
where~$I(a)$ is defined as~$I(a) = \{i \in \{1,\ldots,n\} \, \vert \, a_i = \underline{a} \}$.
\end{example}

\begin{example}[Tori with property {\rm (CS)}]
\label{ex:cstori}
Let~$L = \psi(T(a)) \subset X$ be a Lagrangian torus in a geometrically bounded symplectic manifold having property {\rm (CS)} as in Definition~\ref{def:CS}. Then~$\ce_L \sim \ce_{T(a)}$. Indeed, take the versal deformation~$b \mapsto \psi(T(a+b))$ and note that, for small enough~$b$, all of its members have property {\rm (CS)}, too. By Proposition~\ref{prop:CS}, the claim follows.  
\end{example}

\begin{example}[Toric fibres]
\label{ex:degtoricfibres}
As in Example~\ref{ex:detoricfibres}, let~$X$ be a toric manifold with moment polytope~$\Delta$ defined by vectors~$v_i$. Suppose that~\eqref{eq:detoricfibres} holds on an open dense subset of~$\Int \Delta$. For every toric fibre~$T(x) = \mu^{-1}(x)$, the map~$v_{T(x)}(b) = T(x+b)$ is a versal deformation by Proposition~\ref{prop:momentmapsvd}. We find
\begin{equation}
	\label{eq:degtoricfibres}
	\ce_{T(x)}(b) \sim d_{\rm IA}(x,\pp \Delta) + \min_{i \in I(x)}\{ \langle b , v_i \rangle \}, 
\end{equation}
where~$\sim$ denotes equality on an open dense subset, see Definition~\ref{def:sim}, and where~$I(x) = \{i \, \vert \, \ell_i(x) = d_{\rm IA}(x, \pp \Delta)\}$. Note that~\eqref{eq:degtoricfibres} holds for all \emph{all} toric fibres in~$X$, even those for which~\eqref{eq:detoricfibres} fails. Example~\ref{ex:degprodtori} is a special case of this example.
\end{example}

In all of the examples known to the author, the displacement energy germ of a Lagrangian torus~$L$ has the form
\begin{equation}
	\label{eq:degtorimin}
	\ce_L(b) \sim a + \min\{\langle b , v_1 \rangle, \ldots, \langle b, v_k \rangle\}, \quad
	a> 0, v_i \in H_1(L) \cong \Z^n. 
\end{equation} 
Now let~$L'$ be another torus whose displacement energy germ has the form~\eqref{eq:degtorimin} with~$a' > 0$ and~$v_1', \ldots v_{k'}' \in \Z^n$. If~$L \cong L'$, then by~\eqref{eq:deginv} we find~$a = a'$,~$k = k'$ and that there is~$\Phi \in \GL(n;\Z)$ such that~$\Phi(v_i) = v_i'$ up to relabelling. This explains our interest in the integral geometry of sets of vectors~$v_1,\ldots, v_k \in \Z^n$, see Definition~\ref{def:intindex}.

\begin{remark}
\label{rk:Jholinv}
It would be interesting to interpret~\eqref{eq:degtorimin} in terms of minimal area Maslov two~$J$-holomorphic disks with boundary on~$L$. The quantity~$a > 0$ should correspond to the minimal area and the classes~$v_i \in H_1(L;\Z) \cong \Z^n$ to the boundaries of disks realizing the minimal area~$a$. See~\cite{SheTonVia19} for related ideas.
\end{remark}

\begin{remark}[Virtual displacement energy]
Let~$L$ be a Lagrangian torus whose displacement energy germ has the form~\eqref{eq:degtorimin}. Then~$ve(L) = a > 0$ is an invariant of~$L$ under symplectomorphisms of~$X$, which we call \emph{virtual displacement energy}. In many examples~$ve(L) = e(L)$, but sometimes this is not true. Take for example the monotone torus~$T_{\rm Cl}$ in~$(\CP^2, \omega_{\CP^2})$ with normalization~$\int_{\CP^1} \omega_{\CP^2} = 1$. Then~$e(T_{\rm Cl}) = + \infty$, and
\begin{equation}
	\ce_{T_{\rm Cl}}(b) \sim \frac{1}{3} + \min\left\{ b_1,  b_2 ,  - b_1 - b_2 \right\},
\end{equation}
meaning that~$ve(T_{\rm Cl}) = \frac{1}{3}$. More generally, let~$X$ be toric such that~\eqref{eq:detoricfibres} holds on an open dense subset of~$\Int \Delta$, then~$ve(x) = d_{\IA}(x,\pp \Delta)$ everywhere. Assuming every~$L' \in \cl_L$ has displacement energy germ of the form~\eqref{eq:degtorimin}, we find that~$ve \colon \cl \rightarrow \R$ is continuous, as opposed to~$e(\cdot)$. This is reminiscent of the invariant~$\Psi(\cdot)$ from~\cite{SheTonVia19}. 
\end{remark}

As an application of this discussion and Example~\ref{ex:degtoricfibres}, let us focus on toric manifolds. Note that toric fibres in a toric manifold are the analogue of product tori in~$\C^n$. Therefore, we can call a torus in a toric manifold \emph{exotic} if it is not equivalent to a toric fibre. 

\begin{corollary}
\label{cor:toricexotic}
Let~$X$ be a monotone toric manifold satisfying\footnote{Recall from the discussion in Example~\ref{ex:detoricfibres}, that this condition conjecturally holds for all monotone toric manifolds and has been checked for~$\dim X \leqslant 18$.}~\eqref{eq:detoricfibres} on an open dense subset of~$\Int \Delta$. Then~$X$ contains infinitely many exotic tori, meaning tori distinct from all toric fibres.
\end{corollary}

\proof 
Let~$T(x) \subset X$ be a toric fibre. We can assume that~$T(x)$ is non-monotone, since the tori embedded in the proof of Theorem~\ref{thm:D} are not montone in~$X$. Recall from Example~\ref{ex:degtoricfibres} that the displacement energy germ of~$T(x)$ is of the form~\eqref{eq:degtoricfibres}. Let~$v_i,v_j$ be a pair of vectors appearing in~\eqref{eq:degtoricfibres}. We claim that~$D_{\rm IA}(v_i,v_j) = 1$, where~$D_{\rm IA}$ is the integral index, see Definition~\ref{def:intindex}. Indeed, by monotonicity, the facets associated to~$v_i$ and~$v_j$ have non-empty intersection and thus the claim follows from the Delzant property of the moment polytope~$\Delta$ of~$X$. For the local exotic tori constructed in Theorem~\ref{thm:D} on the other hand, the displacement energy germ is given by~\eqref{eq:degthetaprodalt}, where the minimum is taken over a set containing a pair of vectors for which~$D_{\rm IA} \neq 1$ whenever the corresponding Markov triple~$\fm$ is non-trivial. 
\proofend

\section{Almost toric fibrations}
\label{sec:atfs}

\subsection{Overview}

Almost toric fibrations (ATFs) are a type of Lagrangian torus fibrations on four-dimensional symplectic manifolds introduced by Symington~\cite{Sym03}, which are slightly more general than toric fibrations. We give a brief outline of ATFs here with an emphasis on the properties used in the context of this paper. For more details on almost toric fibrations, we recommend~\cite{Eva23}. 

Toric fibrations lie at the intersection of Hamiltonian group actions and Lagrangian torus fibrations (or integrable systems, the two can be used interchangeably for the present purposes). Indeed the moment map~$\mu \colon X \rightarrow \R^n$ of a toric action on a symplectic manifold~$(X,\omega)$ can be viewed as a special type of Lagrangian torus fibration. It is special in the sense that the moment map yields \emph{global} action coordinates, i.e.\ action coordinates on the regular locus~$\mu^{-1}(\Int \Delta)$ of~$\mu$. In the vicinity of a compact regular fibre, every Lagrangian torus fibration can be equipped with \emph{local} action-angle coordinates compatible with the fibration; this is the classical Arnold--Liouville Theorem~\cite[Theorem 1.42]{Eva23}. However, in general, action-angle coordinates do not extend globally. See Duistermaat~\cite{Dui80} for a study of obstructions to this extension. One obstruction that we encounter when discussing ATFs is the topological monodromy of the torus bundle over the regular locus. If the regular locus is simply-connected, as is the case for toric fibrations, then the topological monodromy is trivial.

To move from toric to almost toric fibrations, we lose the Hamiltonian group action. In terms of Lagrangian torus fibrations/integrable systems, four-dimensional toric fibrations have a very restricted set of allowable singularities, namely only elliptic-regular and elliptic-elliptic singularities appear. We refer to~\cite{Zun96,Zun03} for details. In the almost toric case, the set of permitted singularities is extended by adding the focus-focus case, see~\cite{Zun97}. An almost toric fibration is thus a (proper) Lagrangian torus fibration~$F \colon X^{4} \rightarrow B^2$ with singularities of the three types discussed above. In the context of this paper, we have~$B^2 \subset \R^2$. Note that the focus-focus singularity does not have~$T^2$-symmetry (it only has~$S^1$-symmetry), meaning that we lose the globally defined Hamiltonian~$T^2$-action present in the toric case. Furthermore, the focus-focus singularities map to singular values (called \emph{nodes} or \emph{nodal points}) in the interior of~$B$, meaning that the regular locus~$B_{\rm reg} \subset B$ is not simply connected (as opposed to the toric case) and the torus bundle~$F\vert_{F^{-1}(B_{\rm reg})} \colon F^{-1}(B_{\rm reg}) \rightarrow B_{\rm reg}$ has non-trivial topological monodromy. This means that whenever~$F$ has at least one focus-focus singularity, there are no global action-angle coordinates on~$F^{-1}(B_{\rm reg})$, meaning that a priori it is not obvious what the analogue of the moment map~$\mu$ and the moment polytope~$\Delta$ should be in the almost toric case. However, one can find curves~$C_i$ in~$B$, such that their complement~$B_{\rm reg} \setminus \cup_i C_i \subset B_{\rm reg}$ is simply-connected. In the cases relevant to this paper, every~$C_i$ connects a node to the toric boundary of~$B$. For simplicity, we assume the set of curves~$C_i$ to be mutually disjoint. On~$B_{\rm reg} \setminus \cup_i C_i$, there are global action coordinates. Thus one gets a reasonable analogue of the moment map as follows. First, post-compose~$F\vert_{F^{-1}(B_{\rm reg} \setminus \cup_i C_i)}$ with the diffeomorphism yielding action coordinates. Denote the resulting map by~$\pi_0$. Second, as it turns out, one can choose the curves~$C_i$ in such a way that there is a continuous extension~$\pi$ of~$\pi_0$ mapping~$X$ to the closure~$\Delta_{\pi} = \overline{\im \pi_0} \subset \R^2$ such that~$\Delta_{\pi}$ is a polytope and such that the image of every~$F^{-1}(C_i)$ in~$\Delta_{\pi}$ is a line segment, which we call a \emph{branch cut}. Decorating~$\Delta_{\pi}$ with branch cuts (represented by dashed lines in the conventions of~\cite{Sym03}) indicating the image of the curves~$C_i$ and with nodes (represented by crosses) indicating which fibres contain focus-focus singularities, we obtain a so-called \emph{ATF base diagram} of the Lagrangian torus fibration~$F \colon X \rightarrow B$ we started with. Abusing notation, we denote by~$\Delta_{\pi}$ both the polytope~$\Delta_{\pi} \subset \R^2$ obtained as the image of~$\pi$ and the base diagram, i.e.\ the polytope decorated with branch cuts and nodes. A crucial point in the ATF-framework is that there is a Delzant-type theorem\footnote{The classical Delzant Theorem~\cite{Del88} states that a toric manifold is determined, up to equivariant symplectomorphism, by the integral affine equivalence class of its moment polytope.} for ATF base diagrams. More precisely, the ATF base diagram~$\Delta_{\pi}$ determines the symplectic manifold~$X$ equipped with an almost toric fibration that is unique up to symplectomorphisms which are fibre-preserving up to a union of arbitrarily small neighbourhoods of the nodal points in~$B$, see~\cite[Theorem 8.5]{Eva23}.

The almost toric base diagram~$\Delta_{\pi}$ should be thought of as a geometric-combinatorial object encoding an almost toric fibration~$F \colon X \rightarrow B$ having~$\Delta_{\pi}$ as its ATF base diagram in the sense described above, i.e.\ by removing suitable curves~$C_i \subset B$ and computing action coordinates on their complement. In practice, these subtleties are often glossed over, and, by abuse of terminology, the map~$\pi \colon X \rightarrow \Delta_{\pi}$ is called \emph{almost toric fibration}. When working with the map~$\pi$ explicitly, we keep in mind that it is not smooth over the branch cuts, but that it yields honest action coordinates when restricted to the complement of branch cuts.

Note that, since it involves a choice of curves~$C_i$, the ATF base diagram is \emph{not} uniquely defined up to integral affine transformations. However, the choice of curves~$C_i$ for which the outcome is as described above is very restricted. In fact, there are only two possible choices of~$C_i$ (pointing in opposite directions) for the $i$-th node. Switching at the~$i$-th node from a given choice of~$C_i$ to the one in the opposite direction is called \emph{changing the branch cut}. Note that changing a branch cut does not change the underlying almost toric fibration~$F \colon X \rightarrow B$, merely the way of representing it by a base diagram. Let~$l_i$ be the line containing the branch cut emanating from the~$i$-th node in~$\Delta_{\pi}$. In terms of base diagrams, changing the branch cut at the~$i$-th node corresponds to passing from~$\Delta_{\pi}$ to~$\Delta_{\pi'}$ by the branch cut switching sides in~$l_i$ at the cross representing the~$i$-th node. A computation of action coordinates on the so-obtained simply-connected domain yields that passing from~$\Delta_{\pi}$ to~$\Delta_{\pi'}$ corresponds to applying a piece-wise linear transformation~$\tau \colon \Delta_{\pi} \rightarrow \Delta_{\pi'}$. The map~$\tau$ is continuous, has eigenline~$l_i$ and acts by integral affine transformations on the two halves of~$\R^2 \setminus l_i$. By convention, since everything is up to integral affine transformations, one chooses~$\tau$ to be the identity on one of the halves. Then it is explicitly given by 
\begin{equation}
\label{eq:taumut}
	\tau(x) 
	= 
	\begin{cases}
		\sigma_{w}(x) = x - \det(w \vert x)w, & \det(w \vert x) \geqslant 0 \\
		x, & \det(w \vert x) < 0,
	\end{cases}
\end{equation}
where~$w \in \Z^2$ is the primitive vector contained in~$l_i$ pointing in the direction of the new branch cut. We call~$\tau$ a \emph{mutation of ATF base diagrams}.

Another crucial feature of almost toric fibrations is that, as opposed to toric fibrations, one can apply certain operations to them which produce new ATFs from given ones. These operations are called \emph{nodal trades} and \emph{nodal slides}. Nodal trades exchange an elliptic-elliptic singularity for a focus-focus singularity, merging the two one-parameter families of elliptic-regular singularities emanating from the original elliptic-elliptic singularity. Nodal slides change the position of the node in the ATF base diagram~$\Delta_{\pi}$ by sliding along the line~$l_i$ containing the corresponding node and its branch cut. See~\cite[Sections 8.2/8.3]{Eva23} for details. Thus, given a toric fibration, one can first apply nodal trades to obtain an ATF and then nodal slides to obtain new ATFs on a fixed symplectic manifold. Usually, one changes the branch cut along the way in order for the different branch cuts not to intersect. This process is called \emph{mutation} of ATFs, see~\cite[Section 8.4]{Eva23}. As we shall see, mutations are an effective method, pioneered by Vianna~\cite{Via17, Via16}, to come up with exotic Lagrangian tori.

\subsection{Mutations of the complex projective plane}
\label{ssec:cp2mutations}
For the reader's convenience and to set up notation, let us briefly recall the almost toric fibrations obtained on~$\CP^2$ by applying mutations to the standard toric fibration. For details, we refer to~\cite{Via16} and \cite[Appendix I]{Eva23}. Let~$\CP^2$ be equipped with the symplectic form~$\omega_{\CP^2}$ with normalization~$\int_{\CP^1}\omega_{\CP^1} = \frac{1}{2}$ and the toric system having moment polytope
\begin{equation}	
	\label{eq:cp2momentpolytope}
	\Delta_{\CP^2} = \left\{ x \in \R^2 \, \left\vert \, \langle x,v_i \rangle \geqslant -\frac{1}{3} \right. \right\}, \quad
	v_1 = \begin{pmatrix}
	1 \\ 0
	\end{pmatrix},
	v_2 = \begin{pmatrix}
	0 \\ 1
	\end{pmatrix},
	v_3 = \begin{pmatrix}
	-1 \\ -1
	\end{pmatrix}.
\end{equation}
By applying mutations, one obtains a set of almost toric fibrations which is in bijection with the so-called \emph{Markov triples}, i.e.\ triples~$\fm = (\alpha,\beta,\gamma)$ of natural numbers which satisfy the \emph{Markov equation}
\begin{equation}
	\label{eq:markoveq2}
	\alpha^2 + \beta^2 + \gamma^2 = 3\alpha\beta\gamma.
\end{equation}
These triples appear in a tree rooted in the solution~$(1,1,1)$ (corresponding to the \emph{toric} fibration) which is trivalent except at the first two triples. The edges of the tree correspond to mutations of Markov triples, each of which replaces one of the numbers in the triple by a new Markov number. This comes from the observation that keeping two out of three Markov numbers fixed turns~\eqref{eq:markoveq2} into a quadratic equation having two solutions. On the ATF side, mutating a Markov triple corresponds to mutating the associated almost toric fibration. We denote the polytope of the ATF base diagram corresponding to a Markov triple~$\fm$ by
\begin{equation}
	\label{eq:deltafm}
	\Delta_{\fm} 
	= 
	\left\{ x \in \R^2 \, \left\vert \,  \langle x,u \rangle, \, \langle x,v \rangle, \, \langle x,w \rangle \geqslant - \frac{1}{3} \right. \right\},
\end{equation}
with~$\det(u\vert v) = \alpha^2$,~$\det(v \vert w) = \beta^2$ and~$\det(w\vert u) = \gamma^2$. More precisely, every pair of vectors in~$\{u,v,w\}$ can be mapped by an integral affine transformation to a normal form given by~$(1,0)$ and~$(1-pq,p^2)$, see for example \cite[Lemma I.16]{Eva23}. In this notation,~$p \in \{\alpha,\beta,\gamma\}$ is the corresponding Markov number and~$1 \leqslant q < p$ is a natural number coprime to~$p$. See also~\cite[Figure 2]{Via16} or~\cite[Figure 2]{EvaSmi18}. From this normal form it follows that
\begin{equation}
	\label{eq:markovdifferences}
	u - v = \alpha y_1, \quad
	v - w = \beta y_2, \quad
	w - u = \gamma y_3, 
\end{equation}
where~$y_i \in \Z^2$ are primitive vectors. This will be used in~\S\ref{ssec:viannalift}. For every~$\fm$, there are different ATF base diagrams having the same underlying polytope~$\Delta_{\fm}$, all of which are related by nodal slides. We choose one such base diagram such that all branch cuts are disjoint from the central point~$0 \in \Delta_{\fm}$ and denote the corresponding projection map by~$\pi_{\fm}\colon \CP^2 \rightarrow \Delta_{\fm}$. Recall that this projection map is an honest~$T^2$-moment map on the complement of the branch cuts. Then the almost toric fibre
\begin{equation}
	T_{\fm} = \pi_{\fm}^{-1}(0) \subset \CP^2, 
	\quad
	\fm = (\alpha,\beta,\gamma),
\end{equation}
is a monotone Lagrangian torus called the \emph{Vianna torus} associated to the Markov triple~$\fm$. Vianna~\cite{Via16} proved that one can recover~$\Delta_{\fm}$ up to a~$\GL(2;\Z)$-map from the Maslov two~$J$-holomorphic disks with boundary on~$T_{\fm}$. Since no~$\Delta_{\fm}$ can be mapped to~$\Delta_{\fm'}$ for~$\fm \neq \fm'$ by a~$\GL(2;\Z)$-map (the integral affine angles are~$\GL(2;\Z)$-invariants), one can conclude the following. 

\begin{theorem}[Vianna]
\label{thm:vianna}
The tori~$T_{\fm},T_{\fm'} \subset \CP^2$ are equivalent if and only if~$\fm = \fm'$.
\end{theorem}

This result was conjectured earlier by Galkin--Usnich~\cite{GalUsn10} and proved independently by Galkin--Mikhalkin~\cite{GalMik22}. From their point of view, the tori~$T_{\fm}$ are obtained by considering certain degenerations of~$\CP^2$ to weighted projective spaces~$\CP(\alpha^2,\beta^2,\gamma^2)$. The torus~$T_{\fm}$ can be defined as the preimage under the parallel transport map of the degeneration of the monotone fibre in the weighted projective space. Such degenerations are in bijection with Markov triples by work of Hacking--Prokhorov~\cite{HacPro10}. Note that this viewpoint is closely related to the approach based on ATFs, since one should be able to obtain ATFs by a perturbation of the pull-back to~$\CP^2$ of the orbifold toric system on~$\CP(\alpha^2,\beta^2,\gamma^2)$ under the parallel transport map of the degeneration. 

\subsection{Almost toric fibres}
\label{ssec:atfibres}
In this subsection, we compute the displacement energy of Lagrangian tori appearing as regular fibres of certain almost toric fibrations, which allows us to determine their displacement energy germ. The idea to use the displacement energy germ to distinguish almost toric fibres to obtain results similar to~\cite{Via17, Via16} was suggested to us by Felix Schlenk based on joint work with Yuri Chekanov. Let~$x \in \Delta_{\pi}$ be a point in the complement of the set of nodes in the base diagram of an almost toric fibration~$\pi \colon X \rightarrow \Delta_{\pi}$. Then~$\pi^{-1}(x)$ is a Lagrangian torus. We study the displacement energy of these Lagrangian tori in the special case of almost toric fibrations coming, as for example in~\S\ref{ssec:cp2mutations}, from mutations of an honest toric fibration. The crucial observation is the fact that ATF-fibres which do not lie on a branch cut line are untouched by mutations, see Lemma~\ref{prop:mutationfibreseq}. Hence, the displacement energy of ATF-fibres before a mutation determines the displacement energy of ATF-fibres (away from segments used for nodal slides) after the mutation, see Figure~\ref{fig:3} for an illustration.  

\begin{theorem}
\label{thm:deatfs}
Let~$\pi \colon X \rightarrow \Delta_{\pi}$ be an ATF on a compact monotone symplectic four-manifold, which is obtained by mutations from a toric fibration~$\mu \colon X \rightarrow \Delta$. Then 
\begin{equation}
	e(\pi^{-1}(x)) = d_{\rm IA}(x,\pp \Delta_{\pi})
\end{equation}
for all $x$ in an open dense subset of~$\Delta_{\pi}$.
\end{theorem}

Before turning to the proof of the theorem, we prove two lemmata. 

\begin{lemma}
\label{prop:mutationfibreseq}
Let~$\tau^{-1} \colon \Delta_{\pi} \rightarrow \Delta_{\pi'}$ be a mutation of ATF base diagrams as in~\eqref{eq:taumut} and let~$\pi \colon X \rightarrow \Delta_{\pi}$ be an ATF. Let~$x \in \Delta_{\pi'} \setminus l$ be a point which is not contained in the branch cut line~$l$ of the mutation~$\tau^{-1}$. Then there is an ATF~$\pi' \colon X \rightarrow \Delta_{\pi'}$ such that
\begin{equation}
\label{eq:mutationsfibreseq}
	(\pi')^{-1}(x) = \pi^{-1}(\tau(x)).
\end{equation}
\end{lemma}

\proof 
Nodal slides and nodal trades are supported in arbitrary small neighbourhoods of~$l$ and thus, we can choose~$\pi'$ to differ from~$\pi$ by a change in fibration supported away from~$x \in \Delta_{\pi} \setminus l$ (see for example the proof of~\cite[Theorem 8.10]{Eva23}) and a change of branch cut. Recall that a change of branch cut does not change the underlying fibration, merely its ATF base diagram by applying~$\tau^{-1}$, therefore the claim follows.
\proofend 

\begin{remark}
Actually, a more interesting statement is true. Namely, let~$\pi,\pi'$ be two ATFs related by a nodal trade or nodal slide contained in the branch cut line~$l$, then the almost toric fibres are Hamiltonian isotopic,
\begin{equation}
	\label{eq:mutationsequivfibres}
	(\pi')^{-1}(x) \cong \pi^{-1}(\tau(x)) , \quad
	x \in \Delta_{\pi} \setminus l.
\end{equation}
Here, both almost toric fibrations are fixed a priori, as opposed to the statement in Lemma~\ref{prop:mutationfibreseq}, which comes at the cost of having Hamiltonian isotopies instead of equalities as sets. This statement can be proved by performing toric reduction on segments which are parallel to the branch cut line in question. ATF-fibres project to circles under this reduction and nodal trades and slides correspond to perturbing these reduced fibres, see for example~\cite{Bre22} or~\cite{GroVar21} for related ideas. Area arguments in the two-dimensionsal reduced space allow to deduce~\eqref{eq:mutationsequivfibres}. Although more broadly applicable than Lemma~\ref{prop:mutationfibreseq}, the statement~\eqref{eq:mutationsequivfibres} does not appear in the literature and its proof is more involved, which is why we rely on Lemma~\ref{prop:mutationfibreseq} instead.
\end{remark}

\begin{figure}
  \centering
  \begin{tikzpicture}
    \begin{scope}[shift={(-2,-1.5)}]
      \fill[opacity=0.2] (0,3) rectangle (3,0);
      \draw[very thick] (0,3) -- (0,0) -- (3,0);
      \foreach \x in {0.2,0.4,...,3}
        \draw[ultra thin, opacity=0.6] (\x,3) -- (\x,\x) -- (3,\x);
    \end{scope}
    \draw[->] (1.5,0) -- node[anchor=south] {nodal} node[anchor=north] {trade} (2.5,0);
    \begin{scope}[shift={(3,-1.5)}]
      \fill[opacity=0.2] (0,3) rectangle (3,0);
      \foreach \x in {0.2,0.4,...,3}
        \draw[ultra thin, opacity=0.6] (\x,3) -- (\x,\x) -- (3,\x);
      \draw[white,line width=2pt] (0,0) -- (2,2);
      \draw[very thick] (0,3) -- (0,0) -- (3,0);
      \draw[dashed] (0,0) -- (2,2) node[cross] {};
    \end{scope}
    \draw[->] (6.5,0) -- node[anchor=south] {changing the} node[anchor=north] {branch cut} (8,0);
    \begin{scope}[shift={(8.5,-1.5)}]
      \fill[opacity=0.2] (0,3) rectangle (3,-1);
      \foreach \x in {0.2,0.4,...,3}
        \draw[ultra thin, opacity=0.6] (\x,3) -- (\x,-1);
      \draw[white,line width=2pt] (0,0) -- (2,2);
      \draw[very thick] (0,3) -- (0,-1);
      \draw[dashed] (2,2) node[cross] {} -- (3,3);
    \end{scope}
  \end{tikzpicture}
 
  \caption{A nodal trade at a toric vertex followed by a change of branch cut. Level sets of the displacement energy (under the assumption that \eqref{eq:detoricfibres} holds) of toric fibres before and after the ATF-operations. On the white line segment, the displacement energy cannot be determined by the methods of \S\ref{ssec:atfibres}.}
   \label{fig:3}
\end{figure}

\begin{lemma}
\label{lem:diamutationss}
Let~$\tau^{-1} \colon \Delta_{\pi} \rightarrow \Delta_{\pi'}$ be a mutation of ATF base diagrams which come, via mutations, from the Delzant polytope~$\Delta$ of a compact monotone toric manifold. Then for all~$x \in \Delta_{\pi}$,
\begin{equation}
	\label{eq:diamutations}
	d_{\rm IA}(x , \pp \Delta_{\pi})
	= d_{\rm IA}(\tau^{-1}(x), \tau^{-1} (\pp \Delta_{\pi}))
	=  d_{\rm IA}(\tau^{-1}(x), \pp \Delta_{\pi'}).
\end{equation}
\end{lemma}

\proof
The main idea of the proof is a subdivision into chambers, illustrated in Figure~\ref{fig:4}. Let~$\Delta_{\pi} = \{ x \in \R^2 \, \vert \, \langle x, v_i \rangle +  \lambda  \geqslant 0 \}$ and~$\Delta_{\pi'} = \{ x \in \R^2 \, \vert \, \langle x, v_i' \rangle +  \lambda  \geqslant 0 \}$ with the same conventions as in~\eqref{eq:momentpolytope}. The defining vectors~$v_i$ are labelled starting with the edge adjacent to the branch cut line~$l$ with respect to which the mutation is performed and in the mathematically positive sense. The labelling conventions for the~$v_i'$ is determined below. We denote the edge corresponding to~$v_i$ by~$E_i \subset \Delta_{\pi}$ and the one corresponding to~$v_i'$ by~$E_i' \subset \Delta_{\pi'}$. The fact that the constant~$\lambda > 0$ is the same in every defining equation of~$\Delta_{\pi}$ and~$\Delta_{\pi'}$ follows from monotonicity. Here we suppose that the branch cut line~$l$ of the mutation~$\tau^{-1}$ intersects the interior of an edge opposite to the vertex of the branch cut at which the mutation is performed. In that case the number of edges of~$\Delta_{\pi}$ and of~$\Delta_{\pi'}$ coincides and we denote it by~$N$. The case where~$l$ intersects a vertex is similar up to relabelling some indices. For the mutation~$\tau^{-1} \colon \Delta_{\pi} \rightarrow \Delta_{\pi'}$, where~$\tau$ is given by~\eqref{eq:taumut}, we find that there is~$i_0 \in \{2,\ldots, N-1\}$ and a labelling of the defining vectors~$v_i'$ of~$\Delta_{\pi'}$ such that 
\begin{equation}
	\label{eq:definingvectmut}
	v_1' = \sigma_w^* v_1 \. , \ldots , \.
	v_{i_0}' = \sigma_w^* v_{i_0}, \.\.
	v_{i_0 + 1}' = v_{i_0} \., \ldots , \.
	v_N' = v_{N-1}.
\end{equation}
Furthermore, since the two edges~$E_1,E_N \subset \Delta_{\pi}$ adjacent to the vertex at the branch cut at which the mutation is performed map to the same edge~$E_1' \subset \Delta_{\pi'}$ under~$\tau^{-1}$, we obtain~$v_1' = v_N$. 
For~$i \in \{1,\ldots,N\}$ let
\begin{equation}
	\label{eq:chambers}
	K_i = \{ x \in \Delta_{\pi} \, \vert \, 
	d_{\rm IA}(x, \pp \Delta_{\pi}) 
	= d_{\rm IA}(x,E_i)	
	= \langle x, v_i \rangle + \lambda
	\}
\end{equation}
be the chamber of~$\Delta_{\pi}$ in which the distance to the~$i$-th facet is minimal. The chambers~$K_i' \subset \Delta_{\pi'}$ are defined similarly. We claim that the chamber~$K_i$ is the triangle with vertices~$0,V_i,V_{i+1}$, where~$V_i,V_{i+i}$ are the vertices adjacent to the edge~$E_i$, see Figure~\ref{fig:4}. Indeed, since $d_{\rm IA}(x, \pp \Delta_{\pi}) = \min_i \left\{ d_{\rm IA}(x,E_i) \right\}$ the subset~$ \pp K_i \cap \Int \Delta_{\pi}$ of the boundary of~$K_i$ is formed by the sets satisfying~$d_{\rm IA}(x,E_i) = d_{\rm IA}(x,E_j)$, where~$j$ is the index of an adjacent edge. This is equivalent to~$\langle x, v_i - v_j \rangle = 0$, and the latter equation defines a line containing~$0$ and a vertex of~$\Delta_{\pi}$ adjacent to~$E_i$. The corresponding statements hold for~$\Delta_{\pi'}$. We have
\begin{eqnarray*}
	\tau^{-1}(K_1) \cup \tau^{-1}(K_N) &=& K_1', \\
	\tau^{-1}(K_{i_0}) &=& K_{i_0}' \cup K_{i_0 + 1}', \\
	\tau^{-1}(K_i) &=& K_i', \quad 1 < i < i_0,\\
	\tau^{-1}(K_i) &=& K_{i+1}' \quad i_0 < i < N,
\end{eqnarray*}
which can be deduced from~$\tau^{-1}(0)=0$ together with the fact that~$\tau^{-1}$ induces a bijection of vertices of~$\Delta_{\pi}$ to vertices of~$\Delta_{\pi'}$, except for the vertex~$V_1 \subset \Delta_{\pi}$, which maps to the interior of~$E_1'$ under~$\tau^{-1}$, and the vertex~$V_{i_0}'$ whose preimage lies in the interior of~$E_{i_0}$. Let~$x \in K_i$ with~$1 < i < i_0$. Then 
\begin{eqnarray*}
	d_{\rm IA} (x,\pp \Delta_{\pi})
	&=& \langle x , v_i \rangle + \lambda \\
	&=& \langle x , (\sigma_w^{-1})^*\sigma_w^* v_i \rangle + \lambda \\
	&=& \langle \sigma_w^{-1} x , \sigma_w^* v_i \rangle + \lambda \\
	&=& d_{\rm IA} (\sigma_w^{-1}x, \pp \Delta_{\pi'}) \\
	&=& d_{\rm IA} (\tau^{-1}(x), \pp \Delta_{\pi'}).
\end{eqnarray*}
We have used the definition~\eqref{eq:chambers} of the chambers~$K_i$, equation~\eqref{eq:definingvectmut}, and the fact that $\tau^{-1}(K_i) = K_i'$. For~$x \in K_i$ in another chamber, the claim follows similarly. Since~$\Delta_{\pi} = \cup_i K_i$, the claim of the lemma follows.
\proofend

\proofof{Theorem~\ref{thm:deatfs}}
Let~$S \subset \Delta_{\pi}$ be the union of branch cut lines used in the mutations~$\tau_k, \ldots, \tau_1$ bringing~$\Delta_{\pi}$ to the toric base~$\Delta$, i.e.\ for which we have~$\Delta = (\tau_k \circ \ldots \circ \tau_1)(\Delta_{\pi})$. The set~$S$ is a finite union of line segments. Since~$\mu \colon X \rightarrow \Delta$ is a monotone toric four-manifold,~\eqref{eq:detoricfibres} holds on an open dense subset~$U \subset \Delta$, see Example~\ref{ex:detoricfibres}. Let~$x$ be an element in the open dense subset~$\tau^{-1}(U) \setminus S \subset \Delta_{\pi}$, where~$\tau = \tau_k \circ \ldots \circ \tau_1$. By Lemma~\ref{prop:mutationfibreseq}, the mutations of polytopes can be realized by ATFs such that
\begin{equation}
	\label{eq:atfibresequiv}
	\pi^{-1}(x) = \mu^{-1}(\tau(x)).
\end{equation}
By Lemma~\ref{lem:diamutationss}, we find
\begin{equation}
	e(\pi^{-1}(x))
	=
	d_{\rm IA} (\tau(x), \pp \Delta)
	=
	d_{\rm IA} (x , \tau^{-1}(\pp \Delta))
	=
	d_{\rm IA} (x , \pp \Delta_{\pi}).
\end{equation}
\proofend

\begin{figure}
    \centering
    \begin{subfigure}{0.5\textwidth}
        \centering
        \begin{tikzpicture}[scale=0.45]	
        	\fill[black!15] (-5,0)--(5,0)--(3,4)--(-2,4)--(-5,0);
        	\draw[thick,black!50] (-8,0)--(8,0);
        	\draw[very thick] (-5,0)--(-2,-4)--(3,-5)--(6,-2)--(3,4)--(-2,4)--(-5,0);
        	\draw[very thick,dotted] (0,0)--(-5,0);
        	\draw[very thick,dotted] (0,0)--(-2,-4);
        	\draw[very thick,dotted] (0,0)--(3,-5);
        	\draw[very thick,dotted] (0,0)--(6,-2);
        	\draw[very thick,dotted] (0,0)--(3,4);
        	\draw[very thick,dotted] (0,0)--(-2,4);
        	\node at (7,0.55){$\ell$};
        	\node at (-2,-1.3){$K_1$};
        	\node at (0.2,-2.5){$K_2$};
        	\node at (3,-2.3){$K_3$};
        	\node at (3,0.2){$K_4$};
        	\node at (0.3,2.5){$K_5$};
        	\node at (-2.3,1.6){$K_6$};
        	\node at (5,3){$\Delta_{\pi}$};
        \end{tikzpicture}
    \end{subfigure}%
    ~ 
    \begin{subfigure}{0.5\textwidth}
        \centering
         \begin{tikzpicture}[scale=0.45]	
        	\fill[black!15] (-5,0)--(5,0)--(-3,4)--(-8,4)--(-5,0);
        	\draw[thick,black!50] (-8,0)--(8,0);
        	\draw[very thick] (-8,4)--(-2,-4)--(3,-5)--(6,-2)--(5,0)--(-3,4)--(-8,4);
        	\draw[very thick,dotted] (0,0)--(-8,4);
        	\draw[very thick,dotted] (0,0)--(-2,-4);
        	\draw[very thick,dotted] (0,0)--(3,-5);
        	\draw[very thick,dotted] (0,0)--(6,-2);
        	\draw[very thick,dotted] (0,0)--(5,0);
        	\draw[very thick,dotted] (0,0)--(-3,4);
        	\node at (7,0.55){$\ell$};
        	\node at (-2.8,0){$K_1'$};
        	\node at (0.2,-2.5){$K_2'$};
        	\node at (3,-2.3){$K_3'$};
        	\node at (4.2,-0.75){$K_4'$};
        	\node at (1,1){$K_5'$};
        	\node at (-3.5,2.85){$K_6'$};
        	\node at (2,3){$\Delta_{\pi'}$};
        \end{tikzpicture}
    \end{subfigure}
    \caption{Chambers (delimited by dotted lines) before and after a mutation with eigenline $\ell$. For readablity, we have omitted branch cuts and nodes.}
    \label{fig:4}
\end{figure}

Let~$\Delta_{\pi} = \{\ell_i(x) \geqslant 0\}$, where the functionals~$\ell_i(x) = \langle x, v_i \rangle + \lambda$ are chosen by the same conventions as in~\eqref{eq:momentpolytope}. Then we deduce the following.

\begin{corollary}
\label{cor:degatfs}
Let~$\Delta_{\pi} = \{\ell_i(x) \geqslant 0\}$ be the base of an ATF~$\pi \colon X \rightarrow \Delta_{\pi}$ satisfying the same hypotheses as in Theorem~\ref{thm:deatfs}. Then 
	\begin{equation}
		\label{eq:degatfs}
		\ce_{\pi^{-1}(x)}(b) 
		\sim
		d_{\rm IA}(x, \pp \Delta_{\pi}) + \min_{i \in I(x)}\{\langle b , v_i \rangle\}, \quad
		I(x) = \{i \, \vert \, \ell_i(x) = d_{\rm IA}(x, \pp \Delta_{\pi})\},
	\end{equation}
for all regular fibres~$\pi^{-1}(x) \subset X$. 
\end{corollary}

\proof
Let~$\pi^{-1}(x)$ be a regular fibre. We can assume, possibly after performing changes in branch cuts, that~$x$ does not lie on a branch cut\footnote{For monotone polytopes, one can always make sure that this can be done without producing intersections of branch cuts.}. Recall that the map~$\pi \colon X \rightarrow \Delta_{\pi}$ yields honest action coordinates on the complement of branch cuts. Therefore, for every neighbourhood~$U \subset \Delta_{\pi}$ of~$x$ contained in the complement of branch cuts,~$\pi$ is a~$T^2$-moment map and we obtain a versal deformation~$v_{\pi^{-1}(x)}(b) = \pi^{-1}(x + b)$. This follows from Proposition~\ref{prop:momentmapsvd}. The claim now follows from Theorem~\ref{thm:deatfs}.
\proofend

This implies Vianna's Theorem~\ref{thm:vianna}. Indeed, for the monotone Vianna tori~$T_{\fm} \subset \CP^2$, the displacement energy germ~\eqref{eq:degatfs} recovers the integral affine type of the full polytope~$\Delta_{\fm}$ and thus Theorem~\ref{thm:vianna} follows. 

\begin{example}[Exotic tori in~$\CP^2$]
\label{ex:nonmonvianna}
Furthermore, certain non-monotone almost toric fibres can be distinguished by the same method. Indeed, let~$x$ be a point such that~$d_{\rm IA}(x , \pp \Delta)$ is realized by the integral affine distance to two different edges. Then its displacement energy germ recovers the integral affine geometry of the vertex formed by these two edges. The determinant of the primitive vectors associated to the edges forming this vertex is, by the discussion surrounding~\eqref{eq:deltafm}, equal to the square of the associated Markov number. We conclude that~$\CP^2$ contains infinitely many distinct one-parameter families of non-monotone tori. These families are in bijection with the set of Markov numbers and a triple of such families intersects at a monotone Vianna torus~$T_{\fm}$ if the three Markov numbers form a Markov triple. This implies some of the results in~\cite[Section 7]{SheTonVia19}.
\end{example}

\begin{remark}
In~\cite{Via17}, the results in~$\CP^2$ were extended to other monotone toric four-manifolds. These also follow from the methods discussed here. Although an analogue of Theorem~\ref{thm:deatfs} does not in general hold for ATFs obtained by mutation from \emph{non-monotone} toric manifolds, the displacement energy of almost toric fibres can still be determined by induction on mutations. This way, one can find exotic tori in non-monotone toric four-manifolds. A detailed treatment of this will appear in \cite{Sch24}.
\end{remark}

\begin{remark}
\label{rk:viannaareaeq}
For every Markov triple~$\fm$ the torus~$T_{\fm}$ is area-equivalent to the monotone toric fibre~$T_{\rm Cl} = T_{(1,1,1)}$, meaning that every pair of Vianna tori yields a positive answer to Question~\ref{q:main}. Indeed, as on an open dense subset of its domain, the tori in the versal deformation of~$T_{\fm}$ are Hamiltonian isotopic to toric fibres, we can pick a Lagrangian isotopy from~$T_{\fm}$ to a toric fibre which can be concatenated with a Lagrangian isotopy from that toric fibre to~$T_{\rm Cl}$. Since both endpoints of the resulting Lagrangian isotopy are monotone, we find~$T_{\fm} \sim T_{\rm Cl}$. See also Remark~\ref{rk:areaequivalent}, where the same idea is applied to the tori~$\Upsilon_k \subset \C^3$. 
\end{remark}

\section{Constructions and proofs}
\label{sec:constrproofs}

\subsection{A family of symplectic reductions}
\label{ssec:famsympred}
In this subsection we prepare the construction of the~$\Upsilon_k$-tori given in~\S\ref{ssec:torivd}.
The~$T^3$-action on~$\C^3$ by coordinate-wise rotation is Hamiltonian and generated by the standard moment map,
\begin{equation}
	\mu \colon \C^3 \rightarrow \R_{\geqslant 0}^3, \quad
	\mu(z_1,z_2,z_3) = ( \pi \vert z_1 \vert^2 , \pi \vert z_2 \vert^2 , \pi \vert z_3 \vert^2)
\end{equation}
Let~$k \in \Z_{\geqslant 2}$ and define the moment map
\begin{equation}
	\label{eq:nuk}
	\nu_k \colon \C^3 \rightarrow \Delta_k,\quad
	\nu_k = (\mu_1 + k \mu_3 , \mu_2 - \mu_3).
\end{equation}
This moment map generates a Hamiltonian~$T^2$-action induced by the standard~$T^3$-action via the inclusion
\begin{equation}
	\label{eq:Tk}
	T^2 \hookrightarrow T^3, \quad \exp(\theta_1,\theta_2) \mapsto \exp(\theta_1,\theta_2,k \theta_1 - \theta_2).
\end{equation}
Here,~$\Delta_k = \{ (x,y) \in \R^2 \, \vert \, x \geqslant 0, x+ky \geqslant 0\}$ denotes the image of~$\nu_k$. The set of critical points is given by~$\Crit \nu_k = \cup_{i \neq j}\{z_i = z_j = 0\}$ and the critical values are given by~$\pp \Delta_k \cup (\R_{\geqslant 0} \times \{0\})$. Denote the subtorus obtained as image of the inclusion~\eqref{eq:Tk} by~$T_k \subset T^3$. 

We discuss symplectic reduction with respect to~$\nu_k$ for~$c = (c_1,c_2) \in \Int \Delta_k$, paying special attention to the critical values of the form~$c = (c_1,0)$. Every level set~$\nu_k^{-1}(c)$ contains one orbit with stabilizer isomorphic to~$\Z/k\Z$, namely the orbit obtained by intersecting~$\nu_k^{-1}(c)$ with~$\{z_1 = 0\}$. This orbit maps to an orbifold point of order~$k$ in the reduced space~$\Sigma_{k,c}$. For~$c_2 \neq 0$, the group~$T_k$ acts freely on the complement in~$\nu_k^{-1}(c)$ of that orbit. For~$c = (c_1,0)$, the level set~$\nu_k^{-1}(c)$ intersects~$\{z_2 = z_3 = 0\}$ in the orbit~$\mu^{-1}(c_1,0,0)$ and thus for~$c_1 > 0$, the stabilizer of that orbit is isomorphic to~$S^1$. Removing this exceptional orbit, we can perform symplectic reduction to obtain an orbifold reduced space with one puncture. To summarize, for~$c \in \Int \Delta_k$ we have the following reduction diagram, 
\begin{equation}
\label{eq:rednuk}
\begin{tikzcd}
	Z_{k,c} \arrow[r, hook] \arrow[d, "\pi_{k,c}"] &
	\C^3 \\´
	(\Sigma_{k,c}, \omega_{k,c}),
\end{tikzcd}
\end{equation}
where~$Z_{k,c} = \nu_k^{-1}(c)$ in case~$c_2 \neq 0$, and~$Z_{k,c} = \nu_{k}^{-1}(c) \setminus \mu^{-1}(c_1,0,0)$ in case~$c_2 = 0$ and where~$(\Sigma_{k,c}, \omega_{k,c})$ is the corresponding reduced space. 

\begin{proposition}
\label{prop:redspaces}
For~$c_2 \neq 0$, the reduced space~$\Sigma_{k,c}$ is a symplectic orbifold sphere with one orbifold point of order~$k$. For~$c_2 = 0$, the reduced space is an orbifold sphere with one orbifold point of order~$k$ and one puncture. The quotient space carries an induced effective Hamiltonian~$S^1$-action generated by a moment map~$H_{k,c} \colon \Sigma_{k,c} \rightarrow [0,\alpha_k(c)]$, where 
\begin{equation}
\label{eq:redspacesarea}
	\alpha_k(c)
	=
	\int_{\Sigma_{k,c}} \omega_{k,c}
	=
	\begin{cases}
		\frac{c_1}{k}, & c_2 \geqslant 0 \\
		\frac{c_1}{k} + c_2, & c_2 < 0
	\end{cases}
\end{equation}
is the symplectic area of the reduced space. Furthermore, product tori in~$\C^3$ map to orbits of the~$H_{k,c}$-action on~$\Sigma_{k,c}$. More precisely, we have
\begin{equation}
	\label{eq:orbitlift}
	\pi_{k,c}^{-1}(H_{k,c}^{-1}(c_3)) 
	= 
	\begin{cases}
		T(c_1 - kc_3, c_3 + c_2,c_3), & c_2 \geqslant 0 \\
		T(c_1 - kc_3 + kc_2, c_3 , c_3 - c_2), & c_2 < 0.
	\end{cases}
\end{equation}
\end{proposition}

\proof 
The topology of the reduced spaces has been discussed above. Away from the orbifold point and the puncture in the case~$c_2=0$, the reduction~\eqref{eq:rednuk} is a special case of toric reduction (see the discussion in~\S\ref{ssec:sympred}), meaning that there is an induced Hamiltonian circle action on~$\Sigma_{k,c}$ which is toric. This can be used to prove~\eqref{eq:redspacesarea} and~\eqref{eq:orbitlift}, see Proposition~\ref{prop:redtoricfibres}. 

For the reader's convenience, let us nonetheless discuss the details. The Hamiltonian~$T^3$-action on~$\C^3$ descends to the quotients~$\Sigma_{k,c}$ since~$\nu_k$ (generating a subaction of the standard~$T^3$-action) commutes with~$\mu$, see Proposition~\ref{prop:redgroupaction}. The induced~$T^3$-action on~$\Sigma_{k,c}$ has stabilizer~$T_k$, meaning that it yields an effective Hamiltonian action by~$S^1 \cong T^3 / T_k$ on the quotient. A moment map for this induced action can be defined by picking a circle which is complementary to~$T_k \subset T^3$, e.g.\ the circle~$\{1\} \times \{1\} \times S^1 \subset T^3$, whose action is generated by~$\mu_3$. For $c_2 \geqslant 0$, the moment map~$\mu_3 \vert_{Z_{k,c}} \colon Z_{k,c} \rightarrow \R$ descends to the Hamiltonian~$H_{k,c}\colon \Sigma_{k,c} \rightarrow \R$. In case $c_2 < 0$, we take $\mu_3 \vert_{Z_{k,c}} + c_2$, instead. By computing the image~$\mu_3(Z_{k,c})$, we find that~$H_{k,c}$ takes values in~$[0,\alpha_k(c)]$ (respectively~$(0,\alpha_k(c)]$ for~$c = (c_1,0)$) where~$\alpha_k(c)$ is defined as in~\eqref{eq:redspacesarea}. Since this Hamiltonian circle action is effective, it is toric and thus (as an action coordinate) its range computes the area of~$\Sigma_{k,c}$, proving~\eqref{eq:redspacesarea}. 

To prove~\eqref{eq:orbitlift}, note that~$H^{-1}_{k,c}(c_3)$ is an orbit of the residual action induced by the~$T^3$-action on~$\C^3$. By~$T^3$-equivariance of~$\pi_{k,c}$, see Proposition~\ref{prop:redgroupaction}, this orbit lifts to an orbit of the~$T^3$-action on~$\C^3$, i.e.\ a product torus~$T(a_1,a_2,a_3) \subset \C^3$. To determine~$a_1,a_2,a_3 \in \R_{>0}$, recall that on~$Z_{k,c}$, we have~$\nu_k = c$ and thus~$a_1 + ka_3 = c_1$ and~$a_2 - a_3 = c_2$ by definition of~$\nu_k$, see~\eqref{eq:nuk}. Furthermore,~$c_3 > 0$ is equal to the area of the sublevel set~$H_{k,c}^{-1}(\leqslant c_3)$, since~$H_{k,c}$ generates a Hamiltonian circle action and by its normalization, we find that~$c_3$ is the integral affine distance of the point~$(a_1,a_2,a_3)$ to the union of facets~$\{x_1 = 0\} \cup \{x_2 = 0\}$. This yields~\eqref{eq:orbitlift}.
\proofend

\subsection{Tori and their versal deformations}
\label{ssec:torivd}

We use the same set-up and notation as in~\S\ref{ssec:famsympred}.
 
\begin{definition}
Let~$a_1,a_2 > 0$ with~$a_2 < \alpha_k(a_1,0)$ and set
\begin{equation}
	\Upsilon_k(a_1,a_2) = \pi_{k,(a_1,0)}^{-1}(C(a_2)),
\end{equation}
where~$C(a_2) \subset \Sigma_{k,(a_1,0)}$ is a closed embedded curve in~$\Sigma_{k,(a_1,0)}$ bounding a disk (not containing the puncture) of symplectic area~$a_2$ in the punctured orbifold~$\Sigma_{k,(a_1,0)}$. 
\end{definition}

See also Figure~\ref{fig:1b}. Note that~$C(a_2) \subset \Sigma_{k,(a_1,0)}$ is unique up to Hamiltonian isotopy and thus so is~$\Upsilon_k(a_1,a_2)$, since Hamiltonian isotopies lift by Proposition~\ref{prop:redlifting}. 

\begin{remark}
The same construction can also be carried out, with small modifications at various points, for~$k \in \{-1,0,1\}$. Note that in the cases~$k \in \{-1,0\}$, the reduced space is non-compact and in the case~$k = 1$ it does not have an orbifold point. To simplify the exposition, we do not treat these cases in detail. Furthermore, the tori we obtain in these cases are already known, namely,~$\Theta_{-1} \cong T^3_{\rm Ch}$ is the~$3$-dimensional version of the Chekanov torus,~$\Theta_{0} \cong S^1 \times T^2_{\rm Ch}$ is the product of the Chekanov torus with a circle and~$\Theta_1$ is Hamiltonian isotopic to a monotone product torus and thus not exotic. 
\end{remark}

Let~$a_1,a_2 > 0$ with~$a_2 < \alpha_k(a_1,0)$. From now on, we denote~$\Upsilon_k = \Upsilon_k(a_1,a_2)$ to simplify notation. To construct a versal deformation of~$\Upsilon_k$, we choose a~$T^2$-equivariant Weinstein chart, $\varphi \colon T^*T^3 \supset V \hookrightarrow \C^3$, which intertwines rotation of the first two components of~$T^3$ in~$T^*T^3$ with the~$T_k$-action on~$\C^3$. Rotation of the first two coordinates of~$T^*T^3 = (\R/\Z)^3 \times \R^3 = \{(\theta,p) = (\theta_1,\theta_2,\theta_3,p_1,p_2,p_3)\}$ equipped with~$\omega_{T^*T^3} = \sum_i d\theta_i \wedge dp_i$ is generated by the moment map~$\eta(\theta,p) = (p_1,p_2)$. Since~$\varphi$ is equivariant, we find that $\varphi^* \nu_k = \eta + (a_1,0)$ on~$U$. As in~\eqref{eq:toriweinsteinvd}, we choose the versal deformation
\begin{equation}
	\label{eq:upsilonkvd}
	v_{\Upsilon_k}^{\varphi} \colon H^1(T^3;\R) \cong \R^3 \supset U \rightarrow \cl_{\Upsilon_k}, \quad
	v_{\Upsilon_k}^{\varphi} (b) = \varphi (T^3 \times \{b\}),
\end{equation}
where we have identified~$H^1(T^3,\R) \cong \R^3$ with the fibre part of~$T^*T^3 = T^3 \times \R^3$ and~$U$ is a neighbourhood of~$0 \in \R^3$ such that~$T^3 \times U$ is contained in the domain~$V$ of~$\varphi$.

\begin{proposition}
\label{prop:vdlift}
Let~$b \in U$, then 
	\begin{equation}
		v_{\Upsilon_k}^{\varphi} (b) = \pi_{k,(a_1 + b_1,b_2)}^{-1}(C),
	\end{equation}
for some closed embedded curve~$C \subset \Sigma_{k,(a_1 + b_1,b_2)}$, bounding a disk of area~$a_2 + b_3$.
\end{proposition}

\proof
By equivariance of~$\varphi$, the torus~$v_{\Upsilon_k}^{\varphi}(b) = \varphi (T^3 \times \{b\})$ is~$T_k$-invariant. Furthermore, since~$\varphi^* \nu_k = \eta + (a_1,0)$, we find that~$v^{\varphi}_{\Upsilon_k}(b)$ is contained in the level set~$\nu_k^{-1}(a_1 + b_1, b_2)$. Therefore,~$v^{\varphi}_{\Upsilon_k}(b)$ projects to some compact embedded Lagrangian~$C \subset \Sigma_{k,(a_1 + b_1,b_2)}$. Let us now prove that~$C$ bounds a disk of area~$a_2 + b_3$. Set 
\begin{equation}
	\phi_b = (\varphi \vert_{T^3 \times \{b\}})_* \colon 
	H_1(T^3) \rightarrow H_1(v^{\varphi}_{\Upsilon_k}(b)),
\end{equation}
and let~$e_1,e_2,e_3 \in H_1(T^3)$ be the standard basis of the first homology of~$T^3 = \R^3/\Z^3$. The classes~$\phi_b e_1, \phi_b e_2$ are generated by orbits of the~$T_k$-action and hence~$(\pi_{k,(a_1 + b_1,b_2)}\vert_{v_{\Upsilon_k}^{\varphi}(b)})_* \phi_b e_i = 0$ for~$i = 1,2$. Furthermore, we find that
\begin{equation}
	\left(\pi_{k,(a_1 + b_1,b_2)}\vert_{v_{\Upsilon_k}^{\varphi}(b)} \right)_* \phi_b e_3 = [C] \in H_1(C), 
\end{equation}
where~$[C]$ is the fundamental class with appropriate orientation. Indeed, any basis of $H_1(v_{\Upsilon_k}^{\varphi}(b))$ maps to a basis of~$H_1(C)$ under the reduction map. Recall that~$C$ bounds a disk in~$\Sigma_{k,(a_1+b_1,b_2)}$. Let~$d_b \in H_2(\Sigma_{k,(a_1 + b_1,b_2)},C)$ be the class of this disk with~$\pp d_b = [C]$. Since~$C$ is contractible, the disk~$d_b$ has a unique lift~$D_b \in H_2(\nu_k^{-1}(a_1 + b_1, b_2), v^{\varphi}_{\Upsilon_k}(b))$ with boundary~$\pp D_b = \phi_b e_3$. The cylinder~$C_{e_3}$ swept out by~$\phi_0 e_3$ under a Lagrangian isotopy~$i_t$ (contained in the image of~$\varphi$) with endpoints~$v_{\Upsilon_k}^{\varphi} (0)$ and~$v_{\Upsilon_k}^{\varphi} (b)$ has symplectic area~$\langle [\omega], C_{e_3} \rangle= b_3$. Indeed, this follows from the fact that~$v_{\Upsilon_k}^{\varphi}$ is a versal deformation, see Definition~\ref{def:vd}, together with~\eqref{eq:lagfluxalt}. Furthermore~$\pp C_{e_3} = \phi_b e_3 - \phi_0 e_3$ and therefore it can be compactified to a sphere~$S = D_0 + C_{e_3} - D_b$, by attaching two disks at its boundary. Using the property~\eqref{eq:reductionomega} of symplectic reduction, we obtain 
\begin{equation}
	\langle [\omega] , D_b \rangle
	= \langle [\pi_{k,(a_1 + b_1,b_2)}^* \omega_{k,(a_1 +b_1, b_2)}] , D_b \rangle
	= \langle [\omega_{k,(a_1 +b_1, b_2)}] , d_b \rangle.
\end{equation}
For~$b = 0$, i.e. for the torus~$v_{\Upsilon_k}^{\varphi}(0) = \Upsilon_k$, this quantity is equal to~$a_2$, by the definition of~$\Upsilon_k = \Upsilon_k(a_1,a_2)$. Since~$S \in H_2(\C^3) = 0$, we find 
\begin{equation}
	0 = \langle [\omega] , S \rangle
	= \langle [\omega], D_0 \rangle + \langle [\omega], C_{e_3} \rangle - \langle [\omega], D_b \rangle
	= a_2 + b_3 - \langle [\omega_{k,(a_1 +b_1, b_2)}] , d_b \rangle,
\end{equation}
which proves the claim. 
\proofend

In other words, all tori in the versal deformation~$v_{\Upsilon_k}^{\varphi}$ are lifts of curves in the respective reduced spaces of the action generated by~$\nu_k$ and furthermore, the parameters~$b_1,b_2$ determine the reduced space and~$b_3$ the area of the curve.\smallskip

Note that whenever~$b_2 \neq 0$, the reduced space~$\Sigma_{k,(a_1 + b_1, b_2)}$ is not punctured and hence any two curves~$C,C' \subset \Sigma_{k,(a_1 + b_1, b_2)}$ enclosing the same symplectic area are Hamiltonian isotopic. In particular, the corresponding versal deformations are Hamiltonian isotopic to certain product tori. This is a typical example of the \emph{instability of exotic tori} discussed in~\S\ref{ssec:disttori}. More precisely, combining Propositions~\ref{prop:vdlift} and~\ref{prop:redspaces}, find the following.

\begin{corollary}
\label{cor:prodisom}
For all~$b \in U$ with~$b_2 \neq 0$, the tori in the versal deformation of~$\Upsilon_k$ are Hamiltonian isotopic to the product tori,
\begin{equation}
	\label{eq:prodisom}
	v_{\Upsilon_k}^{\varphi}(b) 
	\cong
	\begin{cases}
		T(a_1 - ka_2 + b_1 - kb_3, a_2 + b_2 + b_3, a_2 + b_3), & b_2 > 0 \\
		T(a_1 - ka_2 + b_1 + kb_2 - kb_3, a_2 + b_3, a_2 - b_2 + b_3), & b_2 < 0. 
	\end{cases} 
\end{equation}
\end{corollary}

\proof
Let~$C \subset \Sigma_{k,(a_1 + b_1, b_2)}$ be the projection of~$v_{\Upsilon_k}^{\varphi}(b)$ as in Proposition~\ref{prop:vdlift}. For~$b_2 \neq 0$, the reduced space is a orbifold sphere with one orbifold point of order~$k$ by Proposition~\ref{prop:redspaces}. Since~$C$ bounds a disk of area~$a_2 + b_3$, it is Hamiltonian isotopic to the orbit~$H_{k,(a_1 + b_1, b_2)}^{-1}(a_2 + b_3)$ of the residual~$S^1$-action on the reduced space. By Proposition~\ref{prop:redlifting}, the Hamiltonian isotopy can be lifted to yield,
\begin{equation}
	v_{\Upsilon_k}^{\varphi}(b) \cong \pi^{-1}_{k,(a_1 + b_1,b_2)} (H_{k,(a_1 + b_1, b_2)}^{-1}(a_2 + b_3)).
\end{equation}
The claim follows from applying~\eqref{eq:orbitlift}.
\proofend

\begin{remark}
\label{rk:supphamisot}
Note that the support of the Hamiltonian isotopy in Corollary~\ref{cor:prodisom} mapping~$v_{\Upsilon_k}^{\varphi}(b)$ to the corresponding product torus can be chosen to lie arbitrarily close to the set~$\nu^{-1}_k(a_1+b_1,b_2)$ on which we perform symplectic reduction. In particular, since~$k \geqslant 2$, its support lies in the ball~$B(a_1+b_1+b_2)$. This will be used in \S\ref{ssec:proofthmB}.
\end{remark}

\subsection{Proof of Theorem~\ref{thm:A}}
\label{ssec:proofthmA}
We compute the displacement energy germ of~$\Upsilon_k = \Upsilon_k(a_1,a_2)$ on the subset~$\{b_2 \neq 0\} \subset H^1(T^3;\R) \cong \R^3$. Note that this is the subset of~$H^1(T^3;\R)$ corresponding to tori which are Hamiltonian isotopic to product tori by Corollary~\ref{cor:prodisom}. This means that we do not compute the displacement energy of the tori~$\Upsilon_k(a_1,a_2)$ (corresponding to versal deformations with~$b_2 = 0$) themselves, see Remark~\ref{rk:deupsilonk}. By~\eqref{eq:prodisom} and~\eqref{eq:deprodtori}, we find 
\begin{equation}
	\label{eq:devdups}
	e(v_{\Upsilon_k}^{\varphi}(b))
	=
	\begin{cases}
		\min\{ a_1 - ka_2 + b_1 - kb_3, a_2 + b_2 + b_3, a_2 + b_3 \}´, & b_2 > 0, \\
		\min\{ a_1 - ka_2 + b_1 + kb_2 - kb_3 , a_2 + b_3, a_2 - b_2 + b_3\}, & b_2 < 0
	\end{cases}
\end{equation}
for~$b \in U$. The germ at the origin of~$b \mapsto e(v^{\varphi}_{\Upsilon_k}(b))$ depends on which terms in the minima in~\eqref{eq:devdups} are dominant for~$b=0$. In case~$a_1 - ka_2 < a_2$, the first terms in the respective minima are dominant. For~$a_1 - ka_2 = a_2$, all terms are dominant. Note however that~$b_3 < b_3 - b_2$ for~$b_2 < 0$ and~$b_3 < b_3 + b_2$ for~$b_2 > 0$, meaning that there are only three dominant terms in total. By the same token, there is only one dominant term in the case~$a_1 - ka_2 > a_2$. We find
\begin{equation}
	\label{eq:upsdegerm}
	\mathcal{E}_{\Upsilon_k}(b)
	=
	\begin{cases}
		a_1 - ka_2 + \min\{\langle b , v_1 \rangle , \langle b , v_2 \rangle \}, & a_1 - ka_2 < a_2 \\
		a_1 - ka_2 + \min\{\langle b , v_1 \rangle, \langle b , v_2 \rangle, \langle b , v_3 \rangle\}, & a_1 - ka_2 = a_2 \\
		a_2 + b_3, & a_1 - ka_2 > a_2,
	\end{cases}
\end{equation}
with 
\begin{equation}
	v_1 =
	\begin{pmatrix}
		1 \\
		k \\
		-k
	\end{pmatrix}, 
	v_2 =
	\begin{pmatrix}
		1 \\
		0 \\
		-k
	\end{pmatrix},
	v_3 =
	\begin{pmatrix}
		0 \\
		0 \\
		1
	\end{pmatrix},
\end{equation}
for all~$b \in \R^3$ with~$b_2 \neq 0$. Let~$a =(a_1,a_2) \in U_k$ and~$a' = (a_1',a_2') \in U_{k'}$, with~$U_k,U_{k'}$ defined as in~\eqref{eq:uk}. This means we are only interested in the first two cases in~\eqref{eq:upsdegerm}. Indeed, note that the displacement energy germ in the case~$a_1 - ka_2 > a_2$ does not depend on~$k$ and thus cannot distinguish the corresponding~$\Upsilon_k$. If~$\Upsilon_k(a) \cong \Upsilon_{k'}(a')$, then~$\ce_{\Upsilon_k} \sim \ce_{\Upsilon_{k'}}$ and, since~\eqref{eq:upsdegerm} is of the form~\eqref{eq:degtorimin}, there is a~$\Phi \in \GL(3;\Z)$ mapping~$\{v_1,v_2\}$ to~$\{v_1',v_2'\}$ in the non-monotone and~$\{v_1,v_2,v_3\}$ to~$\{v_1',v_2',v_3'\}$ in the monotone case. However since~$D_{\rm IA}(v_1,v_2) = k$ by Proposition~\ref{prop:intindex}, this implies~$k = k'$. This proves Theorem~\ref{thm:A}. 
\proofend

\begin{remark}
\label{rk:areaequivalent}
It follows from the proof of Theorem~\ref{thm:A} that~$\Upsilon_k(a_1,a_2) \sim T(a_1 -ka_2, a_2,a_2)$, i.e.\ that every~$\Upsilon_k$ is area-equivalent in the sense of Definition~\ref{def:weakequiv}. Again, we use the fact that small versal deformations of~$\Upsilon_k$ are Hamiltonian isotopic to product tori. The desired Lagrangian isotopy is obtained by concatenation of a deformation of~$\Upsilon_k(a_1,a_2)$  to a suitable torus in its versal deformation with its Hamiltonian isotopy to a product torus and finally with a deformation of the product torus to~$T(a_1 -ka_2, a_2,a_2)$. Such a Lagrangian isotopy has the desired endpoints and vanishing Lagrangian flux. By Remark~\ref{rk:fluxareaclasses}, this proves the claim. 

\end{remark}

\begin{remark}
\label{rk:deupsilonk}
Recall that~$e(\cdot)$ is an upper semicontinuous function on the space of Lagrangians~$\cl$ equipped with the~$C^1$-topology. Since we have computed the displacement energy on an open dense subset of a neighbourhood of~$\Upsilon_k \in \cl$ in the proof of Theorem~\ref{thm:A} above, we deduce~$e(\Upsilon_k(a_1,a_2)) \geqslant a_1 - ka_2$ in case~$a_1 - ka_2 \leqslant a_2$ and~$e(\Upsilon_k(a_1,a_2)) \geqslant a_2$ otherwise. We do not know if these inequalities are sharp. 
\end{remark}

\subsection{Proof of Theorem~\ref{thm:B}} 
\label{ssec:proofthmB}

Let~$(X,\omega)$ be a geometrically bounded symplectic manifold and~$\psi \colon B(R) \rightarrow X$ a symplectic ball embedding. As outlined in~\S\ref{ssec:disttori}, we prove that there is~$\varepsilon > 0$ such that~$\psi(\Upsilon_k)$ has the \emph{locality property}, i.e.\ that~$\ce_{\Upsilon_k} \sim \ce_{\psi(\Upsilon_k)}$ (see Definition~\ref{def:locality}) where~$\Upsilon_k = \Upsilon_k(a_1,a_2)$ for all~$a_1 < \varepsilon$. As in the proof of Theorem~\ref{thm:A} in~\S\ref{ssec:proofthmA}, this allows us to conclude that~$\psi(\Upsilon_k(a)) \cong \psi(\Upsilon_{k'}(a'))$ implies~$k = k'$ whenever~$a \in U_k$ and~$a' \in U_{k'}$. 

Set~$\varepsilon = \min\left\{ \frac{R}{2} , \lambda_S(X,\psi) \right\}$, where~$\lambda_S(X,\psi)$ is defined in~\eqref{eq:lambdapsi}. Let~$a_1 < \varepsilon$ and~$\Upsilon_k = \Upsilon_k(a_1,a_2)$. Note that~$\Upsilon_k \subset B(\varepsilon)$. We apply Lemma~\ref{lem:localitycrit} to the versal deformation~$v_{\Upsilon_k}^{\varphi} \colon U \rightarrow \cl_{\Upsilon_k}$ of~$\Upsilon_k$ constructed in~\S\ref{ssec:torivd}. Here we have restricted, if necessary, the domain~$U$, such that~$a_1 + b_1 + b_2 < \varepsilon$. By Remark~\ref{rk:supphamisot}, this implies that for all~$b \in U$ with~$b_2 \neq 0$,
\begin{equation}
	\label{eq:vdpsiupsilon}
	v_{\Upsilon_k}^{\varphi}(b)
	\cong T(\sigma(b)),
\end{equation}
by a Hamiltonian isotopy supported in~$B(\varepsilon)$. The piecewise affine linear map~$\sigma$ determining~$T(\sigma(b))$ is given as in Corollary~\ref{cor:prodisom}. By Lemma~\ref{lem:localitycrit}, it is sufficient to check that the images~$\psi(T(\sigma(b)))$ of the tori appearing in~\eqref{eq:vdpsiupsilon} have property~{\rm (CS)}, see Definition~\ref{def:CS}. By~\eqref{eq:prodisom}, the condition~$a_1 + b_1 + b_2 < \varepsilon$ implies~$T(\sigma(b)) \subset B(\varepsilon)$. Let~$x = \sigma(b)$. By our choice~$\varepsilon \leqslant \frac{R}{2}$, we find that
\begin{equation}
	x_1 + \ldots + x_n + \underline{x} < \varepsilon + \varepsilon \leqslant \frac{R}{2} + \frac{R}{2} = R
\end{equation} 
and thus~\eqref{eq:CS1} holds. Furthermore,~$\underline{x} < \varepsilon \leqslant \lambda_S(X,\psi)$, and thus~\eqref{eq:CS2} holds. This proves that~$\psi(\Upsilon_k)$ has the locality property and thus that the claim of Theorem~\ref{thm:B} follows.
\proofend

\subsection{Proof of Proposition~\ref{prop:viannalift}}
\label{ssec:viannalift}

Let~$H \colon \C^3 \rightarrow \R$ be the Hamiltonian defined by
\begin{equation}
	H(z_1,z_2,z_3) = \pi (\vert z_1 \vert^2 + \vert z_2 \vert^2 + \vert z_3 \vert^2), 
\end{equation}
which generates the diagonal~$S^1$-action on~$\C^3$. Its symplectic quotient at the level~$c > 0$ is
\begin{equation}
\label{eq:sympredcp2}
\begin{tikzcd}
	H^{-1}(c) = S^5(c) \arrow[r, hook] \arrow[d, "p_c"] &
	\C^3 \\
	(\CP^2, c \omega_{\CP^2}).
\end{tikzcd}
\end{equation}
Let~$T_{\fm} \subset \CP^2$ be the Vianna torus associated with the Markov triple~$\fm$ as in~\S\ref{ssec:cp2mutations} and set
\begin{equation}
	\Theta_{\fm}(a) = p_a^{-1}(T_{\fm}) \subset S^5(a) \subset \C^3. 
\end{equation}

Note that, technically speaking, the identification in~\eqref{eq:sympredcp2} of the reduced spaces~$X_c = H^{-1}(c)/S^1$ with~$\CP^2$ involves a choice of diffeomorphism. Since this will be used in the proof of Lemma~\ref{lem:degliftedvianna}, we make it explicit here. Our definition of~$\CP^2$ is~$\CP^2 = X_1 = H^{-1}(1)/S^1$ equipped with the symplectic form~$\omega_{\CP^2}$ coming from reduction, i.e.\ satisfying~$p_1^*\omega_{\CP^2} = \omega$, where~$\omega$ is the canonical form on~$\C^3$. Note that this coincides with the conventions in~\S\ref{ssec:cp2mutations}, yielding~$\Delta = \Delta_{\CP^2}$ as in~\eqref{eq:cp2momentpolytope} as moment polytope. For~$c_1,c_2 > 0$, consider the natural diffeomorphism~$S_{c_1,c_2} \colon H^{-1}(c_1) \rightarrow H^{-1}(c_2)$ defined by scaling~$z \mapsto \sqrt{\frac{c_2}{c_1}} z$. This induces a diffeomorphism
\begin{equation}
	\label{eq:scalingdiffeo}
	s_{c_1,c_2} \colon X_{c_1} \rightarrow X_{c_2}, \quad
	s_{c_1,c_2}^* \omega_{c_2} = \frac{c_2}{c_1} \, \omega_{c_1}.
\end{equation} 
In~\eqref{eq:sympredcp2}, we have implicitely made use of the diffeomorphism~$s_{1,c} \colon \CP^2 \rightarrow X_c$ with~$s_{1,c}^*\omega_c = c \omega_{\CP^2}$.

\begin{lemma}
\label{lem:degliftedvianna}
\begin{equation}
	\label{eq:degliftvianna}
	\ce_{\Theta_{\fm}(a)} (\ob) \sim d_{\rm IA} (b , (a + b_3)\pp \Delta_{\fm}), \quad
	\ob = (b,b_3). 
\end{equation}
\end{lemma}

Let us first prove that Proposition~\ref{prop:viannalift} follows.\smallskip

\proofof{Proposition \ref{prop:viannalift}}
Let~$\fm = (\alpha,\beta,\gamma)$ be a Markov triple and~$u=(u_1,u_2), v=(v_1,v_2), w=(w_1,w_2) \in \Z^2$ be the defining vectors of the Markov triangle~$\Delta_{\fm}$ as in~\eqref{eq:deltafm}. Lemma~\ref{lem:degliftedvianna} yields,
\begin{equation}
	\label{eq:degliftvianna2}
	\ce_{\Theta_{\fm}(a)}(\overline{b}) 
	\sim
	\frac{a}{3} + \min\{\langle \overline{b}, \overline{u} \rangle, \langle \overline{b} , \overline{v} \rangle, \langle \overline{b}, \overline{w} \rangle\},
\end{equation}
where~$\overline{u} = (u_1,u_2,1)$,~$\overline{v} = (v_1,v_2,1)$ and~$\overline{w} = (w_1,w_2,1)$. We show that by~\eqref{eq:degliftvianna2}, the~$\GL(3;\Z)$-equivalence class of the displacement energy germ~$\ce_{\Theta_{\fm}}$ determines the Markov triple~$\fm = (\alpha,\beta,\gamma)$ which implies the claim. Indeed, the minors of e.g.\ the matrix~$(\overline{u} \vert \overline{v})$ are~$\det(u\vert v) = \alpha^2$,~$u_1 - v_1$ and~$u_2 - v_2$. Therefore Proposition~\ref{prop:intindex} and~\eqref{eq:markovdifferences} allow us to conclude that the integral index is~$D_{\rm IA}(\overline{u}, \overline{v}) = \alpha$. Similarly, we find~$D_{\rm IA} = \beta,\gamma$ for the other two pairs of vectors. 
\proofend

\proofof{Lemma \ref{lem:degliftedvianna}}
To construct a versal deformation of~$\Theta_{\fm}(a) \subset \C^3$, we choose a Weinstein chart~$\varphi \colon T^*T^3 \supset U \rightarrow \C^3$ of~$\Theta_{\fm}(a)$ which is~$S^1$-equivariant with respect to the action by rotation of the last component on~$T^3 = (\R/\Z)^3$ and the action generated by~$H$ on~$\C^3$. Then~$\varphi^*H = p_3 + a$, where~$p_3$ is the third component in~$\R^3$, in the identification~$T^*T^3 = T^3 \times \R^3$. Recall that the associated versal deformation is given by 
\begin{equation}
	v_{\Theta_{\fm}(a)}(\overline{b}) = \varphi\left(T^3 \times \{\overline{b}\}\right), \quad
	\overline{b} = (b,b_3).
\end{equation}
It follows that every torus~$v_{\Theta_{\fm}(a)}(\overline{b}) \subset H^{-1}(a + b_3)$ in the versal deformation is invariant under the action generated by~$H$ and thus projects to a Lagrangian torus~$p_{a+b_3}\left( v_{\Theta_{\fm}(a)}(\overline{b}) \right) \subset X_{a + b_3}$. To compute its displacement energies, we will relate~$v_{\Theta_{\fm}(a)}$ to the versal deformation~$v_{T_{\fm}}(b) = \pi_{\fm}^{-1}(b)$ of~$T_{\fm} \subset X_a \approx \CP^2$ coming from the almost toric fibration~$\pi_{\fm} \colon X_a \rightarrow a\Delta_{\fm}$, discussed (in a more general setting) in~\S\ref{ssec:atfibres}, see in particular the proof of Theorem~\ref{thm:deatfs}. The ATF base diagram is scaled by a factor of~$a > 0$ since we consider the reduced space~$X_a \approx \CP^2$. We will see for example that~$b \mapsto p_a\left( v_{\Theta_{\fm}(a)}(b,0) \right)$ is a versal deformation of~$T_{\fm} = \left( v_{\Theta_{\fm}(a)}(0,0)\right) \subset X_a$. In particular, by the uniqueness of versal deformations up to Hamiltonian isotopy~\eqref{eq:vduniquene}, we find
\begin{equation}
	p_a \left( v_{\Theta_{\fm}(a)}(b,0) \right) \cong v_{T_{\fm}}(\Psi(b)), 
\end{equation}
for some~$\Psi \in \GL(2;\Z)$ and small enough~$b \in \R^2$. More generally, we claim that 
\begin{equation}
	\label{eq:vdinduced}
	w_{b_3}(b) = s_{a +b_3 , a} \left( p_{a + b_3} \left( v_{\Theta_{\fm}(a)} \left( \frac{a + b_3}{a} \, b , b_3 \right) \right) \right) 
	\subset X_{a}
\end{equation} 
is a versal deformation of (a torus equivalent to)~$T_{\fm} \subset X_a$. This is defined for all small enough~$b$. Here~$s_{a + b_3, a} \colon X_{a+ b_3} \rightarrow X_a$ denotes the diffeomorphism~\eqref{eq:scalingdiffeo} satisfying~$s_{a + b_3, a}^* \omega_a = \frac{a}{a+b_3} \omega_{a + b_3}$. As becomes clear below, the scaling by~$\frac{a + b_3}{a}$ in~\eqref{eq:vdinduced} compensates for the scaling of the symplectic form when changing level sets. First, note that~$w_{b_3}(0) \cong T_{\fm} \subset X_a$. Indeed,~$b_3 \mapsto w_{b_3}(0)$ is a Lagrangian isotopy of monotone tori, implying that it has vanishing Lagrangian flux and is therefore induced by a Hamiltonian isotopy. We use~\eqref{eq:lagfluxalt} to prove that~$w_{b_3}$ is a versal deformation in the sense of Definition~\ref{def:vd}. Let~$E_1,E_2,E_3 \in H_1(\Theta_{\fm}(a))$ be the basis induced by the standard basis on~$T^3 = (\R / \Z)^3$ via the Weinstein chart~$\varphi$. Note that~$E_3$ is the class of the~$H$-action and that~$e_1 = \left(p_{a}\vert_{\Theta_{\fm}(a)} \right)_*E_1, e_2 = \left( p_{a}\vert_{\Theta_{\fm}(a)}\right)_*E_2 \in H_1(T_{\fm})$ form a  basis. Since~$v_{\Theta_{\fm}(a)}$ is a versal deformation, we have
\begin{equation}
	\langle \Flux \widehat{\Lambda}_t , E_i \rangle
	= \langle [\omega] , C_{E_i} \rangle
	= b_i,
\end{equation}
for every Lagrangian isotopy~$\widehat{\Lambda}_t$ contained in the image of the versal deformation with~$\widehat{\Lambda}_0 = \Theta_{\fm}(a)$ and~$\widehat{\Lambda}_1 = v_{\Theta_{\fm}(a)}(b)$ and for every relative cycle~$C_{E_i} \in H_2(\C^3, \widehat{\Lambda}_0 \sqcup \widehat{\Lambda}_1)$ contained in the image of~$\varphi$ bounding~$E_i \in H_1(\widehat{\Lambda}_0)$ on~$\Lambda_0$. Now let~$b \in \R^2$ be small enough and let~$\Lambda_t = w_{b_3}(b(t))$ be a Lagrangian isotopy with~$b(0)=0$ and~$b(1)=b$. Since~$(p_{a+b_3})_*E_i = e_i$, every class~$C_{e_i} \in H_2(X_a,\Lambda_0 \sqcup \Lambda_1)$ is the image~$(s_{a+b_3,a} \circ p_{a+b_3})_*C_{E_i}$ of a class
\begin{equation}
	C_{E_i} \in H_2\left(\C^3, v_{\Theta_{\fm}(a)}\left( 0 , b_3 \right) \sqcup v_{\Theta_{\fm}(a)}\left( \frac{a +b_3}{a} b , b_3 \right) \right),
\end{equation} 
contained in~$H^{-1}(a + b_3)$. Since~$v_{\Theta_{\fm}(a)}$ is a versal deformation, we have~$\langle [\omega] , C_{E_i} \rangle = \frac{a + b_3}{a} b_i$ and thus we find,
\begin{eqnarray*}
	\langle \Flux \Lambda_t , e_i \rangle
	&=& \langle [\omega_a] , C_{e_i} \rangle \\
	&=& \langle [p_{a+b_3}^* s_{a+b_3,a}^* \omega_{a}] , C_{E_i} \rangle \\
	&=& \frac{a}{a+b_3} \langle [p_{a+b_3}^*\omega_{a+b_3}] , C_{E_i}\rangle \\
	&=& \frac{a}{a+b_3} \langle [\omega] , C_{E_i}\rangle  \\
	&=& \frac{a}{a+b_3}  \frac{a + b_3}{a} b_i = b_i,
\end{eqnarray*}
where we have used~\eqref{eq:reductionomega} in the fourth equality.
Thus we have proved that~\eqref{eq:vdinduced} is a versal deformation of~$w_{b_3}(0) \cong T_{\fm}$. By the uniqueness of versal deformations~\eqref{eq:vduniquene}, we deduce that there is~$\Psi \in \GL(2;\Z)$ such that 
\begin{equation}
	s_{a +b_3 , a} \left( p_{a + b_3} \left( v_{\Theta_{\fm}(a)} \left( \frac{a + b_3}{a} \, b , b_3 \right) \right) \right) 
	\cong 
	v_{T_{\fm}}(\Psi(b))
\end{equation}
Using~$s_{a +b_3 , a} \circ p_{a+b_3} = p_a \circ S_{a + b_3, a} $, we find
\begin{equation}
	\label{eq:vdtheta0}
	v_{\Theta_{\fm}(a)} \left( \Psi^{-1}(b) , b_3 \right)
	\cong 
	S_{a +b_3 , a}^{-1} \left( p_a^{-1} \left( v_{T_{\fm}} \left( \frac{a}{a+b_3} \, b \right) \right) \right).
\end{equation}
As discussed in~\S\ref{ssec:atfibres}, there is a piece-wise linear map~$\tau = (\tau_1,\tau_2) \colon \R^2 \rightarrow \R^2$, the linear pieces of which are in~$\GL(2;\Z)$, such that
\begin{equation}
	\label{eq:vdtheta1}
	v_{T_{\fm}}(b) \cong \mu^{-1}(\tau(b)), \quad b \in V,
\end{equation}
where~$\mu^{-1}(\cdot)$ denotes toric fibres of the standard moment map~$\mu$ and~$V$ is an open dense subset of a neighbourhood of the origin. See the proof of Theorem~\ref{thm:deatfs} for details. Under~$p_a$, toric fibres lift to product tori in~$\C^3$, 
\begin{equation}
	\label{eq:vdtheta2}
	p_a^{-1}(\mu^{-1}(y_1,y_2))
	=
	T \left( \frac{a}{3} + y_1 , \frac{a}{3} + y_2 , \frac{a}{3} - y_1 - y_2 \right), \quad
	(y_1,y_2) \in a \Delta, 
\end{equation}
see Proposition~\ref{prop:redtoricfibres}. Product tori are mapped to product tori under the scaling map~$S_{a+b_3,a}$, 
\begin{equation}
	\label{eq:vdtheta3}
	S_{a+b_3,a}(T(x)) = T\left( \frac{a}{a+b_3} x \right), \quad
	x \in \R_{\geqslant 0}^3.
\end{equation}
Applying~\eqref{eq:vdtheta1},~\eqref{eq:vdtheta2} and~\eqref{eq:vdtheta3} to~\eqref{eq:vdtheta0}, we find, for every~$b \in V$,
\begin{equation}	
	\label{eq:vdthetafm}
	v_{\Theta_{\fm}(a)} \left( \Psi^{-1}(b) , b_3 \right)
	\cong 
	T\left( \frac{a+b_3}{3} + \tau_1(b), \frac{a+b_3}{3} + \tau_2(b), \frac{a+b_3}{3} - \tau_1(b) - \tau_2(b) \right). 
\end{equation}
Therefore, using Lemma~\ref{lem:diamutationss} we find (recall Definition~\ref{def:sim}) that
\begin{eqnarray*}
	\ce_{\Theta_{\fm}(a)}(\overline{b})
	&\sim & \frac{a + b_3}{3} + \min\{\tau_1(b),\tau_2(b), - \tau_1(b) - \tau_2(b)\} \\
	& = & d_{\rm IA}(\tau(b),(a+b_3)\pp \Delta) \\
	& = & d_{\rm IA}(b,(a+b_3)\tau^{-1}(\pp \Delta)) \\
	& = & d_{\rm IA}(b,(a+b_3)\pp \Delta_{\fm})
\end{eqnarray*}

\proofend

\begin{remark}
\label{rk:3Datf}
One can consider constructing the analogue of an almost toric fibration for~$\C^3 \setminus \{0\}$, starting from a given ATF~$\pi_{\fm} \colon \CP^2 \rightarrow \Delta_{\fm}$, by considering the projection map
\begin{equation}
	\overline{\pi}_{\fm} \colon \C^3 \setminus \{0\} \rightarrow \R^3, \quad
	\overline{\pi}_{\fm}(z) = (h \cdot \pi_{\fm}(p_h(z)), h), \quad
	h = H(z) = \pi (\vert z_1 \vert^2 + \vert z_2 \vert^2 + \vert z_3 \vert^2).
\end{equation}
The image of~$\overline{\pi}_{\fm}$ is the infinite cone over~$\Delta_{\fm}$,
\begin{equation}
	C(\Delta_{\fm}) = \{ (h x_1, h x_2, h) \in \R^3 \, \vert \, h > 0, (x_1,x_2) \in \Delta_{\fm} \}.
\end{equation}
The associated Lagragian torus fibration has three one-parameter families of focus-focus-regular singularities, projecting to three rays~$C(\Delta_{\fm})$. Note that the branch cuts are two dimensional wedges obtained by taking the cones over the three branch cuts in~$\Delta_{\fm}$. On the complement of these, the map~$\overline{\pi}_{\fm}$ generates an honest~$T^3$-action. The torus~$\Theta_{\fm}(a)$ is the fibre over the point~$(0,0,a)$. As in the case of (almost) toric fibres in dimension four, the displacement energy germ of~$\ce_{\Theta_{\fm}(a)}$ is given by the integral affine distance to the boundary of the base of the Lagrangian torus fibration, 
\begin{equation}
	\ce_{\Theta_{\fm}(a)} (\ob) \sim d_{\rm IA} (\ob, \pp C(\Delta_{\fm})).
\end{equation}
See also the discussion in~\S\ref{ssec:atfshigherdim}. 
\end{remark}

\subsection{Proof of Theorem~\ref{thm:C}} 
\label{ssec:proofthmC}

Let~$(X^6,\omega)$ be a geometrically bounded symplectic manifold and~$\psi \colon B(R) \rightarrow X$ a Darboux chart. Let~$\varepsilon = \min\left\{ \frac{3R}{4} , \lambda_S(X,\psi) \right\}$. As in the proof of Theorem~\ref{thm:B} in~\S\ref{ssec:proofthmB}, we show that embeddings~$\psi(\Theta_{\fm})$ of lifted Vianna tori~$\Theta_{\fm} = \Theta_{\fm}(a)$ with~$a < \varepsilon$ have the locality property,~$\ce_{\psi(\Theta_{\fm})} \sim \ce_{\Theta_{\fm}}$. We apply Lemma~\ref{lem:localitycrit} to the versal deformation~$v_{\Theta_{\fm}} \colon V \rightarrow \cl_{\Theta_{\fm}}$ constructed in~\S\ref{ssec:viannalift}. Since~$\Theta_{\fm} \subset B(\varepsilon)$, we restrict the domain~$V$ such that the Lagrangian tori in its image contained in~$B(\varepsilon)$. Recall from the proof of Proposition~\ref{prop:viannalift} in~\S\ref{ssec:viannalift} that~$v_{\Theta_{\fm}}(\overline{b}) \cong T(\sigma(\overline{b}))$ for~$\overline{b}$ in an open dense subset of~$V$. Since this Hamiltonian isotopy is lifted from the corresponding reduced space~$H^{-1}(c)/S^1 \cong (\CP^2, c \omega_{\CP^2})$ with~$c < \varepsilon$, it has support in~$B(\varepsilon)$ and thus it suffices to check that the image~$\psi(T(\sigma(x)))$ of the tori appearing the versal deformation have property {\rm (CS)} to apply Lemma~\ref{lem:localitycrit}. Since~$\psi(T(\sigma(x))) \subset B(\varepsilon)$, this follows from our choice of~$\varepsilon$. See also the end of the proof of Theorem~\ref{thm:B} in~\S\ref{ssec:proofthmB}. 

It remains to show that~$\psi(\Theta_{\fm}(a)) \sim \psi(\Theta_{\fm'}(a))$ for all Markov triples~$\fm, \fm'$, meaning that these provide an answer to Question~\ref{q:main}. Indeed, recall from Remark~\ref{rk:viannaareaeq} that~$T_{\fm} \sim T_{(1,1,1)}$, which implies, by lifting the corresponding Lagrangian isotopy to~$\C^3$, that~$\Theta_{\fm}(a) \sim T\left( \frac{a}{3} , \frac{a}{3}, \frac{a}{3} \right) \subset \C^3 $ for every Markov triple~$\fm$. Since this Lagrangian isotopy is supported in~$H^{-1}(a)$, it can be carried out in the image of~$\psi$ and the claim follows. \proofend

\subsection{Proof of Theorem~\ref{thm:D}}
\label{ssec:proofthmD} In the same vein as above, we first construct exotic tori in~$\C^n$ and then prove that they have the locality property with respect to the displacement energy germ.
Let~$n > 3$. For a Markov triple~$\fm$ and~$a, a_4, \ldots, a_n > 0$, let,
\begin{equation}
	\label{eq:thetaprod}
	\Theta_{\fm}^n
	=
	\Theta_{\fm}^n(a,a_4,\ldots,a_n)
	=
	\Theta_{\fm}(a) \times T(a_4,\ldots,a_n) \subset \C^n.
\end{equation}

\begin{lemma}
\label{lem:degthetaprod}
For~$\frac{a}{3}\leqslant  \min\{a_4,\ldots,a_n\}$, the displacement energy germ of~$\Theta_{\fm}^n$ is given by
\begin{equation}
	\label{eq:degthetaprod}
	\ce_{\Theta_{\fm}^n}(B)
	\sim
	\min_{i \in I} \left\{ d_{\rm IA}(\overline{b},\pp\Delta_{\fm}) , \frac{a}{3} + b_i  \right\}, \quad
	I  = \left\{ i \in \{4,\ldots,n\} \, \left\vert \, a_i = \frac{a}{3} \right.\right\}, 
\end{equation}
where~$B = (b_1,\ldots,b_n) = (\overline{b},b_4,\ldots,b_n) \in \R^n$.
\end{lemma}

As in the proof of Proposition~\ref{prop:viannalift}, see in particular~\eqref{eq:degliftvianna2}, this can be written as
\begin{equation}
	\label{eq:degthetaprodalt}
	\ce_{\Theta_{\fm}^n}(B)
	\sim
	\frac{a}{3} +
	\min_{i \in I}\{\langle B , U \rangle, \langle B, V \rangle, \langle B, W \rangle, \langle B , e_i \rangle \},
\end{equation}
where~$U = (\overline{u},0), V = (\overline{v},0), W = (\overline{w},0) \in \Z^n$, where~$\overline{u},\overline{v},\overline{w} \in \Z^3$ are defined as in the proof surrounding~\eqref{eq:degliftvianna2}. Again, we find~$D_{\rm IA}(U,V) = \alpha, D_{\rm IA}(U,V) = \beta, D_{\rm IA}(U,V) = \gamma$, meaning that~$\Theta_{\fm}^n(a,a_4,\ldots,a_n) \cong \Theta_{\fm'}^n(a',a_4',\ldots,a_n')$ implies~$\fm = \fm'$ whenever~$\frac{a}{3} \leqslant \min\{a_4,\ldots,a_n\}$ and~$\frac{a'}{3} \leqslant \min\{a_4', \ldots, a_n'\}$. In other words, we find infinitely many non-equivalent tori in~$\C^n$. Recall from~\S\ref{ssec:proofthmC} that~$\Theta_{\fm}(a) \sim T(\frac{a}{3},\frac{a}{3},\frac{a}{3})$, which implies that
\begin{equation}
	\label{eq:thetaprodareaeq}
	\Theta_{\fm}^n(a,a_4,\ldots,a_n) \sim T\left(\frac{a}{3},\frac{a}{3},\frac{a}{3},a_4,\ldots,a_n \right),
\end{equation}
meaning the~$\Theta_{\fm}^n$ provide infinitely many examples answering question~\ref{q:main}. In case~$a_4 = \ldots = a_n = \frac{a}{3}$, the torus~$\Theta_{\fm}^n$ is monotone. 

\begin{remark}
\label{rk:codimoneee}
Note that~\eqref{eq:thetaprod} defines a~$(n-2)$-parametric family of exotic tori. This family lies at the intersection of multiple~$(n-1)$-parametric families of exotic tori. Indeed, recall that the monotone Vianna torus~$T_{\fm}$ lies a the intersection of one-parameter families of exotic tori, see Example~\ref{ex:nonmonvianna} and~\cite[Section 7]{SheTonVia19}, which lift to two-parameter families of exotic tori in~$\C^3$, see Remark~\ref{rk:nonmonotonelifts}. Taking products with product tori, one obtains~$(n-1)$-parameter families of exotic tori in~$\C^n$.
\end{remark}

\proofof{Lemma~\ref{lem:degthetaprod}}
Let~$v_{\Theta_{\fm}(a)}$ be the versal deformation of~$\Theta_{\fm}(a)$ used in the proof of Lemma~\ref{lem:degliftedvianna} in~\S\ref{ssec:viannalift}, where we can assume, up to a change of basis, that~$\psi = \id$. By taking the product with the natural versal deformation of product tori from Example~\ref{ex:degprodtori}, we obtain the following product versal deformation of~$\Theta_{\fm}^n = \Theta^n_{\fm}(a,a_4,\ldots,a_n)$,
\begin{equation}
	\label{eq:vdthetaprod}
	v_{\Theta^n_{\fm}}(B)
	=
	v_{\Theta_{\fm}(a)}(\overline{b}) 
	\times
	T(a_4 + b_4, \ldots, a_n + b_n), \quad
	B = (\overline{b},b_4,\ldots,b_n).
\end{equation}
From~\eqref{eq:vdthetafm}, we deduce, on an open dense subset,
\begin{equation}
	\label{eq:vdthetaprodisom}
	v_{\Theta^n_{\fm}}(B)
	\cong 
	T\left( \frac{a + b_3}{3} + \tau_1(b) , \frac{a + b_3}{3} + \tau_2(b) , \frac{a + b_3}{3} - \tau_1(b) - \tau_2(b) , a_4 + b_4, \ldots, a_n + b_n \right),
\end{equation}
where~$\tau$ is the piecewise integral affine transformation coming from mutations (see~\S\S \ref{ssec:cp2mutations}-\ref{ssec:atfibres}) mapping~$\Delta_{\fm}$ to the toric moment polytope~$\Delta_{(1,1,1)}$. As in the proof of Lemma~\ref{lem:degliftedvianna}, we obtain 
\begin{equation}
	\label{eq:devdthetaprod}
	e(v_{\Theta^n_{\fm}}(B))
	\sim 
	\min\left\{ d_{\rm IA}(\overline{b}, \pp \Delta_{\fm}), a_4 + b_4, \ldots, a_n + b_n \right\}.
\end{equation}
Taking the germ in the case~$\frac{a}{3} \leqslant \min\{a_4,\ldots,a_n\}$ proves the claim. 
\proofend
	
\begin{remark}
\label{rk:asmall}
The condition~$\frac{a}{3} \leqslant \min\{a_4,\ldots,a_n\}$ is necessary for our methods to work, since otherwise, the~$d_{\rm IA}(\cdot, \pp \Delta_{\fm})$ does not appear in the germ of~\eqref{eq:devdthetaprod}. We do not know whether it is a necessary condition for the tori to be inequivalent. 
\end{remark}

\proofof{Theorem \ref{thm:D}}
Let~$(X,\omega)$ be a geometrically bounded symplectic~$2n$-manifold and~$U \subset X$ and open subset. Furthermore, we fix a symplectic ball embedding~$\psi \colon B(R) \hookrightarrow U$. We choose $\Theta_{\fm}^n(a,a_4,\ldots,a_n)$ such that the product tori appearing in its versal deformation have property~{\rm (CS)} with respect to~$\psi$. This allows us to apply Lemma~\ref{lem:localitycrit} and deduce that~$\psi(\Theta_{\fm}^n)$ has the locality property. Let~$a,a_4,\ldots,a_n > 0$ such that~$\frac{a}{3} \leqslant \min\{a_4,\ldots,a_n\}$ and
\begin{equation}
	\label{eq:thetaprodCS}
	\frac{4a}{3} + a_4 + \ldots + a_n < R, \quad
	\frac{a}{3} < \lambda_S(X,\psi).
\end{equation}
The versal deformation~\eqref{eq:vdthetaprod} induces the versal deformation~$v_{\psi(\Theta_{\fm}^n)} = \psi \circ v_{\Theta_{\fm}^n}$. By~\eqref{eq:thetaprodCS}, we can choose the domain~$V$ of this versal deformation small enough such that
\begin{align}
	\frac{4a}{3} + a_4 + \ldots &+ a_n + b_3 + \ldots + b_n \nonumber \\
	\label{eq:vdthetaprodCS1}
 &+ \min_{i \in I}\left\{ \frac{b_3}{3} + \tau_1(b) , \frac{b_3}{3} + \tau_2(b), \frac{b_3}{3} - \tau_1(b) - \tau_2(b) , b_i \right\} < R, \\
	\label{eq:vdthetaprodCS2}
	\frac{a}{3} &+ \min_{i \in I}\left\{ \frac{b_3}{3} + \tau_1(b) , \frac{b_3}{3} + \tau_2(b), \frac{b_3}{3} - \tau_1(b) - \tau_2(b) , b_i \right\} < \lambda_S(X,\psi),
\end{align}
for all~$B = (b_1,\ldots,b_n) = (b,b_3, \ldots, b_n) \in V$.
Since the support of the Hamiltonian isotopy~\eqref{eq:vdthetaprodisom} can be chosen to lie in an arbitrarily small neighbourhood of~$S^5(a+b_3) \times T(a_4 + b_4, \ldots, a_n + b_n)$, the condition~\eqref{eq:vdthetaprodCS1} implies that~$v_{\psi(\Theta_{\fm}^n)}(B) \cong \psi(T(y))$, with support in~$\psi(B(R))$, where~$y$ is given by~\eqref{eq:vdthetaprodisom}. Furthermore~\eqref{eq:vdthetaprodCS1} and~\eqref{eq:vdthetaprodCS2} imply that the versal deformations of~$\Theta_{\fm}^n$ have property {\rm (CS)}. This means that we can apply Lemma~\ref{lem:localitycrit} to conclude that~$\psi(\Theta_{\fm}^n)$ has the locality property. Together with~\eqref{eq:thetaprodareaeq}, this allows us to deduce that if~$\fm \neq \fm'$, then
\begin{equation}
	\psi(\Theta_{\fm}^n (a,a_4,\ldots,a_n)) \ncong
	\psi(\Theta_{\fm'}^n (a,a_4,\ldots,a_n)), \quad
	\psi(\Theta_{\fm}^n (a,a_4,\ldots,a_n)) \sim
	\psi(\Theta_{\fm'}^n (a,a_4,\ldots,a_n)),
\end{equation}
which proves Theorem~\ref{thm:D}.
\proofend

\section{Questions and outlook}
\label{sec:qandp}

\subsection{Questions}
\label{ssec:questions}

The following question came up in discussions with L. Polterovich. 

\begin{question}
\label{q:1}
When are exotic tori unstable? 
\end{question}

As mentioned in~\S\ref{ssec:disttori}, we say that an exotic torus~$L \subset \C^n$ is \emph{unstable} if an open dense subset of tori in its versal deformation are not exotic, meaning they are equivalent to product tori. In a general symplectic manifold~$X$, we can ask if there is a pair of Lagrangians~$L,L' \subset X$ which are \emph{stably distinct}. By this we mean that~$L \ncong L'$ and~$L \sim L'$ such that their versal deformations~$v_L \colon U \rightarrow \cl$ and~$v_{L'} \colon U' \rightarrow \cl$ are distinct in the sense that~$v_{L}(b) \ncong v_{L'}(b')$ for all~$b \in V, b' \in V'$, where~$V \subset U,V' \subset U'$ are neighbourhoods of the origin.

\begin{remark}
\label{rk:stablyexotic}
In \cite{HinZha23}, Hind and Zhang have shown that there are stably exotic tori in balls, ellipsoids and in polydisks in $\C^2$, see \cite[Theorems 2.5-2.7]{HinZha23}. Note however that the authors conjecture that these tori lose their exoticity when considered in all of $\C^2$. 
\end{remark}

All exotic tori in compact ambient spaces known to the author are unstable. We heavily use the instability of the~$\Upsilon_k$- and~$\Theta_{\fm}$-tori to compute their displacement energy germs. Every exotic torus known to the author extends to a family of exotic tori of at least codimension one in its deformation space, see for example Remarks~\ref{rk:nonmonotonelifts} and \ref{rk:codimoneee} and~\S\ref{ssec:higherdim}.

\begin{question}
\label{q:2}
Does every exotic torus come at least in a codimension one family of exotic tori?
\end{question}

Note that the Vianna tori (and their lifts) lie at the intersection of three codimension one families of exotic tori. This suggests a stratified structure of the set of exotic tori in the space of all tori. 

\begin{question}
\label{q:3}
Is the displacement energy germ of every Lagrangian torus of the form~\eqref{eq:degtorimin}?
\end{question}

Again, this is the case in all examples known to the author. One way of answering Question~\ref{q:3} may be to link the displacement energy germ to an invariant coming from Maslov two~$J$-holomorphic disks of minimal area as suggested in Remark~\ref{rk:Jholinv}. Let~$L \subset X$ be a (not necessarily monotone) Lagrangian torus.

\begin{question}
\label{q:4}
Is there a symplectic invariant of~$L$ based on counting Maslov two $J$-holomorphic disks (of minimal area) which has the locality property? Can~$\ce_L$ be recovered from this invariant?
\end{question}

In case~$X$ is monotone and (almost) toric, the displacement energy germ~$\ce_L$ of the monotone fibre~$L$ is of the form~\eqref{eq:degtorimin} and it recovers the geometry of the toric base or of the ATF base diagram, respectively, in the sense that the base polytope is given by
\begin{equation}
	\label{eq:deltacel}
	\Delta_{\ce_L} = \{ x \in \R^n \, \vert \, \langle x , v_i \rangle + a \geqslant 0 \},
\end{equation}
where~$v_1,\ldots,v_n \in \Z^n$ and~$a > 0$ are given by~\eqref{eq:degtorimin}.

\begin{question}
\label{q:5}
Let~$L \subset X$ be a monotone Lagrangian torus with displacement energy germ~$\ce_L$ given by~\eqref{eq:degtorimin}. Is there a Lagrangian torus fibration with base~$\Delta_{\ce_L}$ as defined in~\eqref{eq:deltacel}?
\end{question}

Note that in that case~$\ce_L$ is the germ of the function~$d_{\rm IA}(\cdot,\pp \Delta_{\ce_L})$. In that sense, one might expect a correspondence between the set of monotone Lagrangian tori in~$X$ and certain Lagrangian torus fibrations (up to a reasonable equivalence relation, e.g.\ nodal slides in the almost toric case) on~$X$. In Remark~\ref{rk:3Datf} and~\S\ref{ssec:atfshigherdim} we discuss some evidence for a positive answer to Question~\ref{q:5} in the case of the~$\Theta_{\fm}$- and~$\Upsilon_k$-tori constructed in this paper. In a second step, one can ask which singularities appear in the Lagrangian torus fibration and how these are related to properties of~$L \subset X$.

\subsection{Dimension four}
\label{ssec:dimensionfour}
In dimension four there is, as of now, no analogue of the local exotic tori constructed in this paper. It is still open whether every Lagrangian torus~$L \subset \C^2$ is equivalent to either a product or a Chekanov torus, see~\cite{Riz19} for evidence in that direction. In case this classification holds, there can be no such analogue due to a lack of exotic tori, even in~$\C^2$. Furthermore, embeddings of the Chekanov torus by Darboux charts are not always exotic in a robust sense. The following example makes use of symmetric probes, see~\cite{Bre23}, and will appear as an application of the classification of toric fibres in non-monotone~$S^2 \times S^2$ in~\cite{BreKim23}.

\begin{example}
\label{ex:cheknonexotic}
Let $\lambda > 2$ and let~$X_{\lambda}$ be the symplectic manifold obtained by equipping~$S^2 \times S^2$ with the symplectic form~$(\lambda\omega_{S^2}) \oplus \omega_{S^2}$. Recall that~$T_{\rm Ch}(a) \subset \C^2$ is defined as the lift of a curve enclosing area~$a>0$ in the reduced space fibering over the diagonal in the toric moment map image of~$\C^2$, see~\S\ref{ssec:introkey}. Hence, whenever $a < r$, we find a copy $T_{\rm Ch}(a) \subset B^4(r)$. In \cite{BreKim23}, it is shown that for every ball embedding $\varphi \colon B^4(r) \hookrightarrow X_{\lambda}$, the image $\varphi(T_{\rm Ch}(a))$ is Hamiltonian isotopic to a \emph{split torus} in $S^2 \times S^2$, meaning a torus splitting as a product of circles. Thus Chekanov tori in non-monotone $S^2 \times S^2$ are not exotic. This stands in contrast with the results of this paper, see Corollary \ref{cor:toricexotic}. Furthermore, again in \cite{BreKim23}, it is shown that for most (but not all) $a$, the embedded torus $\varphi(T_{\rm Ch}(a))$ is Hamiltonian isotopic to $\varphi(T(a,b))$ for some non-monotone product torus $T(a,b) \subset \C^2$. By choosing~$a > 0$ small enough, one can arrange for~$\varphi(T_{\rm Ch}(a))$ to have property {\rm (CS)}. Again, this example stands in contrast to the local exotic tori in dimension six and higher constructed in this paper, see Remark~\ref{rk:cstori}.
\end{example}

However, in a recent work joint with Johannes Hauber and Joel Schmitz~\cite{BreHauSch23}, we treat a \emph{semi-local} analogue of the ideas developed in this paper yielding tori which are exotic in a more robust way. Instead of embedding Lagrangian tori contained in symplectic balls, we treat Lagrangian tori contained in more complicated symplectic domains. For example, let~$L \subset X^4$ be a Lagrangian sphere in a symplectic four-manifold. Then a neighbourhood of~$L$ is symplectomorphic to a neighbourhood of the zero-section in~$T^*S^2$ which serves as our domain. Any neighbourhood of the zero-section in~$T^*S^2$ contains three types of tori, namely product, Chekanov and Polterovich tori. For a discussion of the Polterovich torus, we refer to~\cite{AlbFra08}, see also \cite{OakUsh16} for a detailed discussion in higher dimensions. Under similar smallness assumptions on the tori as in Theorems~\ref{thm:B} and~\ref{thm:C}, see also Remark~\ref{rk:epsilon}, one can prove that these types of tori are not equivalent in~$X^4$. More generally, we consider embeddings of certain Liouville domains~$B_{d,p,q}$ into~$X^4$, yielding~$d+1$ distinct tori. For example~$B_{2,1,1} \cong T^*S^2$ and~$B_{1,2,1} \cong T^*\RP^2$, see~\cite{Eva19} for a discussion of the spaces~$B_{d,p,q}$. These spaces appear as Milnor fibres of cyclic quotient singularities, or equivalently, neighbourhoods of \emph{Lagrangian pinwheels}, see~\cite{EvaSmi18} for more details. The space~$B_{d,p,q}$ admits an almost toric fibration with~$d$ nodes lying on one branch cut line and the tori studied in~\cite{BreHauSch23} appear as regular fibres over the branch cut line of these ATFs. For a fixed area class~$a > 0$, the corresponding torus changes equivalence class whenever a node is moved across it. As in the present paper, in~\cite{BreHauSch23} we prove a locality property which implies that the tori remain distinct after embedding them into a geometrically bounded~$X^4$. 

The main difference with the purely local case of Darboux embeddings is that there are constraints from symplectic topology for symplectic embeddings of the type~$B_{d,p,q} \supset U \rightarrow X^4$ to exist. A full discussion of this can be found in~\cite{EvaSmi18} for the case~$d=1$ for~$X = \CP^2$.

\subsection{Almost toric fibrations in higher dimensions}
\label{ssec:atfshigherdim}

In the same vein as Question~\ref{q:5}, one can try to construct a Lagrangian torus fibration on~$\C^3$ having the~$\Upsilon_k$-tori as fibres. Note that for the Chekanov torus~$T_{\rm Ch} \subset \C^2$, this works and yields an almost toric fibration on~$\C^2$. We suspect that this leads to a notion of almost toric fibrations in arbitrary dimension. Let us give a rough outline of this.

One can construct a Lagrangian torus fibration by hand, using a global version of the arguments we have used to deal with the versal deformation of the~$\Upsilon_k$. As projection map, consider 
\begin{equation}
	\label{eq:pikG}
	\pi_{k,G} \colon \C^3 \rightarrow \Delta_k \times \R_{\geqslant 0}, \quad
	\pi_{k,G}(z) = (\nu_k(z), \pi_{k,\nu_k(z)}^*G_{k,\nu_k(z)}(z)).
\end{equation}
Since~$\nu_k$ generates the~$T_k$-action considered in~\S\ref{ssec:famsympred}, we can interpret~$\pi_{k,G}$ as a~$T_k$-equivariant perturbation of the toric system on~$\C^3$. Indeed, the last component of~$\pi_{k,G}$ is a lift of a Hamiltonian~$G_{k,c} \colon \Sigma_{k,c} \rightarrow \R$ from the reduced space~$\Sigma_{k,c}$ of the~$T_k$-action. Note that if we lifted~$H_{k,c}$ instead, then we would end up with the standard toric system on~$\C^3$, since the fibres of~$\pi_{k,H}$ are product tori by Proposition~\ref{prop:redspaces}. However, by lifting perturbations of~$H_{k,c}$ having the circles~$C(\cdot) \subset \Sigma_{k,c}$ appearing in the definition of the~$\Upsilon_k$-tori as level sets, we end up with an Lagrangian torus fibration having~$\Upsilon_k$-tori as fibres. A good candidate for a family~$\{G_{k,c}\}_{c \in \Delta_k}$ of functions to define~$\pi_{k,G}$ seems to be: start with~$H_{k,c}$ and perturb it on a small neighbourhood of the ray~$\{(c_1,0)\} \subset \Delta_k$. Just as a similar~$S^1$-equivariant perturbation of the toric system on~$\C^2$ produces an almost toric nodal trade at the vertex of~$\R_{\geqslant 0}^2$, this~$T_k$-equivariant procedure produces an analoguous \emph{nodal trade at the edge}~$\{x_2 = x_3 = 0\} \subset \R^3_{\geqslant 0}$. Singular fibres appear where the level sets of~$G_{k,(c_1,0)}$ intersect the punctures of~$\Sigma_{k,(c_1,0)}$ and they are of focus-focus-regular type, since locally the nodal trade at the vertex is equivalent to starting with the standard toric system on~$\C^2 \times T^*S^1$ and performing a nodal trade in the~$\C^2$-component. 
Nodal slides correspond, in this picture, to varying the level of~$G_{k,(c_1,0)}$ at which the puncture is intersected. 
The corresponding ATF base diagrams are computed by removing two-dimensional branch sets and computing action coordinates on their complement. 

In~$\R_{\geqslant 0}^3$, there are three edges at which one can perform nodal trades and for each of those, one can choose~$k \in \Z_{\geqslant 0}$. Therefore, for each triple~$(k_1,k_2,k_3) \in \Z_{\geqslant 2}^3$, we obtain a series of ATFs and mutations of~$\C^3$. We call~$(k_1,k_2,k_3)$ the \emph{nodal trade vector}. Note that for the nodal trade vector~$(k_1,k_2,k_3) = (2,2,2)$, the~$T^2$-actions at chosen at each edge intersect in a common circle, namely the diagonal circle~$S^1 \subset T^3$. Therefore, the corresponding mutation in~$\C^3$ are compatible with symplectic reduction by the diagonal circle and project to mutations of~$\CP^2$. This reproduces the Markov mutations in~$\CP^2$, see~\S\ref{ssec:cp2mutations}. In other words, our construction for~$k_i = 2$ reproduces the cone construction discussed in Remark~\ref{rk:3Datf} and thus fibres of the higher mutations correspond, in this case, to the tori~$\Theta_{\fm}$. On the other hand, simple mutations at only one edge correspond to the tori~$\Upsilon_k$ as we have seen. Therefore these conjectural ATF-mutations in higher dimensions have fibres which generalize both the~$\Upsilon_k$- and~$\Theta_{\fm}$-tori, with $\Theta_{\fm}$-tori appearing as fibres of ATFs obtained by \emph{higher mutations for the fixed nodal trade vector}~$(k_1,k_2,k_3) = (2,2,2)$ and $\Upsilon_k$-tori appearing as fibres of ATFs obtained by \emph{first order mutations for the nodal trade vectors}~$(k_1,k_2,k_3)= (k,\varnothing,\varnothing)$, where~$\varnothing$ stands for not making a nodal trade. For every~$(k_1,k_2,k_3) \in \Z_{\geqslant 2}^3$, we obtain a Markov-like tree of Lagrangian tori. It would be highly interesting to give a proper construction of these families tori and distinguish them. As in the four-dimensional analogue in~\S\ref{ssec:atfibres}, this should be possible by the use of versal deformations. In particular, it should be possible to find out how the displacement energy germ behaves under mutations and that it can be read off from the base of the corresponding ATF. 

\begin{remark}
Note that this could lead to a theory of ATFs for compact toric manifolds of arbitrary dimension. Indeed, every compact toric manifold~$X$ is, via the Delzant construction, see~\cite{Gui94}, an equivariant symplectic quotient of some~$\C^N$ by a subtorus~$K \subset T^N$ acting on~$\C^N$. We can restrict our attention to nodal trades at codimension~$(k \geqslant 2)$-faces of~$\R^N_{\geqslant 0}$ such that their corresponding tori contain~$K$. Just as the~$(2,2,2)$-ATFs on~$\C^3$ induce ATFs on~$\CP^2$, this general recipe should induce ATFs on~$X$.
\end{remark}

\begin{remark}
Note that the Lagrangian torus fibrations in~\eqref{eq:pikG} are perturbations of the toric moment map which are somewhat \emph{mild} in the sense that, locally, a~$T^2$-action is preserved. What happens if one decides to only preserve an~$S^1 \subset T^3$? This line of inquiry should be started by analzing what perturbations of the toric system on~$\C^2$ make sense in this context, since the reduced spaces are two dimensional in this case.
\end{remark}

\newpage

\bibliographystyle{abbrv}
\bibliography{/home/joe/Documents/Math/000_bibtex/mybibfile}

\end{document}